\documentclass{tac}

\usepackage{amssymb}
\usepackage{amsmath}
\usepackage{array}
\usepackage{bbm}
\usepackage{enumerate}
\usepackage[latin1]{inputenc}
\usepackage[T1]{fontenc}
\usepackage{shadow}
\usepackage[all,cmtip,2cell]{xy}
\UseAllTwocells
\usepackage{graphicx}
\usepackage{graphbox}
\usepackage{caption}
\usepackage{subcaption}
\usepackage{float}
\usepackage{multicol}
\usepackage{tikz}
\usepackage{tikz-cd}
\usepackage{longtable}
\usepackage{color}
\usetikzlibrary{matrix,arrows,decorations.pathmorphing}
\usetikzlibrary{calc}
\usepackage{pdfpages}
\usepackage[colorlinks=true]{hyperref}
\hypersetup{allcolors=[rgb]{0.1,0.1,0.4}}

\tikzset{curve/.style={settings={#1},to path={(\tikztostart)
    .. controls ($(\tikztostart)!\pv{pos}!(\tikztotarget)!\pv{height}!270:(\tikztotarget)$)
    and ($(\tikztostart)!1-\pv{pos}!(\tikztotarget)!\pv{height}!270:(\tikztotarget)$)
    .. (\tikztotarget)\tikztonodes}},
    settings/.code={\tikzset{quiver/.cd,#1}
        \def\pv##1{\pgfkeysvalueof{/tikz/quiver/##1}}},
    quiver/.cd,pos/.initial=0.35,height/.initial=0}

\newcommand{\A}{\mathcal{A}}
\newcommand{\B}{\mathcal{B}}
\newcommand{\CC}{\mathcal{C}}
\newcommand{\QQ}{\mathcal{Q}}
\newcommand{\VV}{\mathcal{V}}
\newcommand{\W}{\mathcal{W}}
\newcommand{\DD}{\mathcal{D}}
\newcommand{\PP}{\mathcal{P}}
\newcommand{\G}{\mathbb{G}}
\newcommand{\N}{\mathbb{N}}
\newcommand{\sd}{D^s}

\newcommand{\one}{\mathbbm{1}}

\mathrmdef{Fun}
\mathrmdef{Hom}
\mathrmdef{gSet}
\mathrmdef{Set}
\mathrmdef{id}
\mathrmdef{Map}
\mathrmdef{Adj}
\mathrmdef{Alg}
\mathrmdef{Comp}
\mathrmdef{Cat}
\mathrmdef{pd}
\mathrmdef{Tree}
\mathrmdef{Cell}
\mathrmdef{I}
\mathrmdef{II}
\mathrmdef{III}
\mathrmdef{IV}
\mathrmdef{V}
\mathrmdef{VI}
\mathrmdef{lax}
\mathrmdef{oplax}
\mathrmdef{Psh}

\definecolor{myblack}{RGB}{0,0,0}
\definecolor{myred}{RGB}{255,0,0}
\definecolor{myblue}{RGB}{0,0,255}
\definecolor{mygreen}{RGB}{0,120,0}
\definecolor{mypurple}{RGB}{120,0,120}
\definecolor{myorange}{RGB}{255,120,0}
\definecolor{mylightblue}{RGB}{0,255,255}
\definecolor{myolive}{RGB}{120,120,0}
\definecolor{myteal}{RGB}{0,120,120}

\newcommand{\red}{\textcolor{myred}{red}}
\newcommand{\blue}{\textcolor{myblue}{blue}}
\newcommand{\green}{\textcolor{mygreen}{green}}
\newcommand{\purple}{\textcolor{mypurple}{purple}}
\newcommand{\orange}{\textcolor{myorange}{orange}}
\newcommand{\lightblue}{\textcolor{mylightblue}{light blue}}
\newcommand{\olive}{\textcolor{myolive}{olive}}
\newcommand{\teal}{\textcolor{myteal}{teal}}

\title{String diagrams for $4$-categories and fibrations of mapping $4$-groupoids}  

\author{Manuel Ara\'{u}jo}     
\address{}
\eaddress{manuel.araujo@tecnico.ulisboa.pt}


\copyrightyear{2024}

\thanks{I am grateful to Pedro Boavida, Christopher Douglas, John Huerta and Roger Picken, Richard Garner and the anonymous referee for many useful comments and suggestions on draft versions of this paper. This work was partially supported by the FCT project grant \textbf{Higher Structures and Applications}, PTDC/MAT-PUR/31089/2017.}

\begin{document}

\maketitle                 

\begin{abstract}
	We introduce a string diagram calculus for strict $4$-categories and use it to prove that given a cofinite inclusion of $4$-categorical presentations, the induced restriction functor on mapping spaces to a fixed target strict $4$-category is a fibration of strict $4$-groupoids. 
\end{abstract}

\section{Introduction}

In this paper we introduce a string diagram calculus for strict $4$-categories and we use it to prove that given an inclusion of $4$-categorical presentations $\PP\hookrightarrow\QQ$, where $\QQ$ is obtained from $\PP$ by adding a finite number of generating cells, the induced restriction functor on mapping spaces to a fixed target strict $4$-category is a fibration of strict $4$-groupoids. The string diagram calculus for strict $4$-categories is based on the string diagram calculus for $n$-sesquicategories introduced in \cite{sesqui_comp} and the \textit{homotopy generators} introduced in \cite{data_struct_quasistrict}. This string diagram calculus and the result on fibrations are essential ingredients in \cite{adj3} and \cite{adj4} on coherence for adjunctions in $3$ and $4$-categories. All results are stated in the context of strict $4$-categories, but the proofs will work without change in the context of a theory of semistrict $4$-categories admitting an appropriate string diagram calculus, which is currently being developed.

\subsection{Results}

Now we state the main result in this paper, for which we need to introduce some terminology. By a \textbf{strict $n$-groupoid} we mean a strict $n$-category all of whose morphisms are \textbf{weakly invertible}. A functor between strict $n$-groupoids is a \textbf{fibration} if one can lift any $k$-morphism along it, extending a given lift of its source. We use the word \textbf{presentation} to mean a $5$-computad for the monad $T_4$ on $4$-globular sets whose algebras are strict $4$-categories. Given a presentation $\PP$ we denote by $F(\PP)$ the $4$-category generated by $\PP$. Given a $4$-category $\CC$, we denote by $\Map(\PP,\CC)$ the $4$-groupoid whose objects are funtors $F(\PP)\to\CC$ and whose $k$-morphisms ($k=1,\cdots, 4$) are the weakly invertible $k$-transfors (also known as natural transformations, modifications and perturbations, for $k=1,2,3$).

\begin{theorem}\label{main}

Let $\CC$ be a strict $4$-category, $\PP$ a presentation and $\QQ$ another presentation, obtained by adding a finite number of cells to $\PP$. Then the restriction map $$\Map(\QQ,\CC)\to\Map(\PP,\CC)$$ is a fibration of strict $4$-groupoids.

\end{theorem} 

\subsection{Related work}

The string diagram calculus for strict $4$-categories builds on the string diagram calculus for $n$-sesquicategories introduced in \cite{sesquicat} and \cite{sesqui_comp} and uses the \textbf{homotopy generators} introduced in \cite{data_struct_quasistrict}. One can therefore implement all the string diagram calculations in this paper in the proof assistant \textbf{Globular} (\cite{globular}).

In \cite{adj3} and \cite{adj4}, we use the results from the present paper to prove \textbf{coherence} results for \textbf{adjunctions} in strict $3$ and $4$-categories. This can then be used to give a simplified proof of the result in \cite{araujo_thesis} on coherence for $3$-dualizable objects in strict symmetric monoidal $3$-categories. This is related via the cobordism hypothesis to the ultimate goal of this project, which to establish a finite presentation for the symmetric monoidal $3$-category of framed cobordisms of dimension $\leq 3$. 

Another approach to string diagrams for higher categories is the theory of associative $n$-categories (\cite{dorn_thesis}), developed into the theory of \textbf{manifold diagrams} in \cite{manifold_diagrams}. The combinatorial counterpart of this theory is equivalent to the zigzags introduced in \cite{high_level_methods}, which forms the basis for the proof assistant \textbf{homotopy.io} (\cite{homotopyio}).

In \cite{folk} the \textbf{folk model strucuture} on the category of strict $n$-categories is introduced. Given $n$-groupoids $E$ and $B$, it is natural to ask whether a map $f:E\to B$ is a fibration in the sense of the present paper if and only if it is a folk fibration. If one could show the two notions of fibration agree, then it would be possible to give a short proof of Theorem \ref{main} for all $n$, using results of \cite{al20}. See Remark \ref{folkfib} for more on this. In the present paper, we give an explicit string diagram proof of Theorem \ref{main}, which has the advantage that it will also apply in any model of semistrict $4$-categories admitting a string diagram calculus. 

Note that in \cite{bg89} the authors construct a model structure on the category of strict $n$-groupoids. However, they define a strict $n$-groupoid as a strict $n$-category where every $k$-morphism has a strict inverse, rather than a weak one. See also \cite{am11}.

\subsection{Future work}

We are working on a theory of \textbf{semistrict $4$-categories} based on string diagrams. Once this is in place, the contents of Sections \ref{functorcat}, \ref{scompfun} and \ref{sfibration} will extend to that setting.
 
\thirdleveltheorems

\section{Background}

We give a brief overview of the string diagram calculus for $n$-sesquicategories, introduced in \cite{sesquicat} and \cite{sesqui_comp}. Then we establish some basic terminology about equivalences in and between strict $n$-categories.

\subsection{Globular pasting diagrams and strict $n$-categories}

Let $\G$ be the globe category, with set of objects $\mathbb{N}$ and morphisms generated by $s,t\colon k\to k+1$ satisfying the usual globularity relations $ts=ss$ and $tt=st$. Let $\gSet:=\Fun(\G^{op},\Set)$ the category of \textbf{globular sets}. Let $\G_n$ be the full subcategory of $\G$ on objects $\{0,\cdots,n\}$ and $\gSet_n:=\Fun(\G_n^{op},\Set)$ the category of \textbf{$n$-globular sets}. 

Recall (e.g. \cite[Chapter 8]{operads_cats}) that one can define a \textbf{strict $n$-category} as an algebra over a monad $T_n:\gSet_n\to\gSet_n$. We now briefly recall the construction of this monad.

\begin{notation}

Denote by $\langle n \rangle$ the totally ordered set  $\{1\leq\cdots\leq n\}$ for $n\in\N$.
 
\end{notation}

\begin{definition} 

A \textbf{globular $k$-pasting diagram} $\pi$ is a diagram $$\xymatrix@1{\langle\ell_k\rangle\ar[r] & \cdots \ar[r] & \langle\ell_1\rangle\ar[r] & \langle 1 \rangle}$$ of totally ordered sets and order preserving maps. Denote by $\pd(k)$ the set of globular $k$-pasting diagrams. The \textbf{source} $s(\pi)$ and the \textbf{target} $t(\pi)$ are both equal to the globular $(k-1)$-pasting diagram obtained by truncation. Thus $\pd$ becomes a globular set.

\end{definition}

Diagrams of totally ordered sets as above can be pictured as trees, with height $k$ and $\ell_i$ nodes at level $i$. They can also be pictured as globular pasting diagrams. We explain this in a couple of examples.

A globular $1$-pasting diagram is just a natural number. The pasting diagram representation of $m$ is a string of $m$ composable arrows $$\xymatrix{\bullet\ar[r] & \bullet\ar[r] & \cdots\ar[r] & \bullet}.$$ 

The tree of height $2$ on the left corresponds to the globular $2$-pasting diagram on the right. \begin{center}\begin{tabular}{lcccr}\includegraphics[scale=2,align=c]{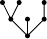} & & & & $\xymatrix@1{\bullet\ruppertwocell\rlowertwocell\ar[r] & \bullet \ar[r] &\bullet \rtwocell & \bullet}$\end{tabular}\end{center}

To each $k$-pasting diagram $\pi$, it is possible to associate a $k$-globular set $\widehat{\pi}\in\gSet_k$ whose cells are the the ones sugested by its pictorial representation.  If we take $\pi$ to be the $2$-dimensional example above, then $\widehat{\pi}$ has four $0$-cells, six $1$-cells and three $2$-cells, and the source and target maps can be easily read from the picture. A precise definition can be found in \cite{operads_cats}. 

\begin{definition}
	
	Given a globular $k$-pasting diagram $\pi$ and a $k$-globular set $X$, we define an \textbf{$X$-labelling} of $\pi$ to be a map of $k$-globular sets $\hat{\pi}\to X$.
	
\end{definition}

This corresponds to the idea of labelling each dot and arrow in the graphical depiction of the pasting diagram, with each label being a cell in $X$ of the correct dimension, satisfying source and target compatibility. 

Consider the functor

$$T_{n}\colon \gSet_n\to\gSet_n$$ defined by $$T_{n}(X)_{k}=\coprod_{\pi\in\pd(k)}\gSet_n(\hat{\pi},X).$$ This functor can be given a monad structure, by interpreting a pasting diagram labelled by pasting diagrams as a pasting diagram.

\begin{definition}
	
	An \textbf{$n$-category} is an algebra over the monad $T_n$. A \textbf{functor} between $n$-categories is a morphism of algebras over $T_n$. 
	
\end{definition} 

This agrees with the standard definitions of strict $n$-categories and strict functors between them.

\subsection{Simple string diagrams and n-sesquicategories}

In \cite{sesquicat}, we defined another monad $$T_n^{\sd}:\gSet_n\to\gSet_n$$ on $n$-globular sets, whose algebras we called \textbf{$n$-sesquicategories}. This is based on a notion of \textbf{simple string diagram}, which plays the same role as globular pasting diagrams in the definition of $T_n$. 

\begin{definition} 

A \textbf{simple $k$-string diagram} $\pi$ is a diagram $$\xymatrix@1{\langle\ell_k\rangle\ar[r] & \cdots \ar[r] & \langle\ell_1\rangle\ar[r] & \langle 1 \rangle}$$ of totally ordered sets and maps which are not required to be order preserving. Denote by $\sd(k)$ the set of simple $k$-string diagrams. The \textbf{source} $s(\pi)$ and the \textbf{target} $t(\pi)$ are both equal to the simple $(k-1)$-string diagram obtained by truncation. Thus $\sd$ becomes a globular set.

\end{definition}

Diagrams of totally ordered sets as above can be pictured as \textit{trees with crossings}, with height $k$ and $\ell_i$ nodes at level $i$. They can also be pictured as \textit{simple string diagrams}, as we now explain. We always read $k$-string diagrams from left to right when $k$ is odd and from top to bottom when $k$ is even.

A simple $1$-string diagram is just a natural number. The string diagram representation of $m$ consists of $m$ dots on a line, which we interpret at $m$ composable morphisms. For example, for $m=2$ we have $$\includegraphics[scale=1.5]{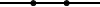}.$$

The tree with crossings on the left corresponds to the simple $2$-string diagram on the right.

\begin{center}\begin{tabular}{lcccr}\includegraphics[scale=2,align=c]{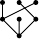} & & & & \includegraphics[scale=1.5,align=c]{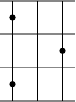}\end{tabular}\end{center} We interpret this diagram as a vertical composition of three $2$-morphisms, each of which consists of a whiskering of a $2$-morphism with two $1$-morphisms.

The tree with crossings on the left corresponds to the simple $3$-string diagram on the right.

\begin{center}\begin{tabular}{lcccr}\includegraphics[scale=2,align=c]{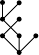} & & & & \includegraphics[scale=1.5,align=c]{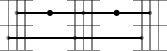}\end{tabular}\end{center}

The tree with crossings on the left corresponds to the simple $4$-string diagram on the right.

\begin{center}\begin{tabular}{lcccr}\includegraphics[scale=2,align=c]{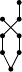} & & & & \includegraphics[scale=1.5,align=c]{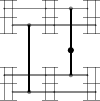}\end{tabular}\end{center}

\begin{definition}

We define a map of globular sets $\pi:\sd\to\pd$ by ``grabbing the leaves and undoing the crossings''. See \cite{sesquicat} for a precise definition.

\end{definition}

\begin{example}

The map $\pi$ on an example simple $2$-string diagram.

\begin{center}\begin{tabular}{ccccc}\includegraphics[scale=2,align=c]{pasting/tree2.pdf} & & $\mapsto$ & & \includegraphics[scale=2,align=c]{pasting/tree.pdf}

\\

\\

\includegraphics[scale=1.5,align=c]{stringdiag/2d_2.pdf}
 & & $\mapsto$ & & $\xymatrix@1{\bullet\ruppertwocell\rlowertwocell\ar[r] & \bullet \ar[r] &\bullet \rtwocell & \bullet}$
               
\end{tabular}\end{center}

\end{example}

\begin{example}
 
We give some more examples of simple string diagrams and their associated globular pasting diagrams.

\begin{center}
 \begin{tabular}{clc}
  \includegraphics[align=c,scale=1.5]{stringdiag/1d.pdf} & $\mapsto$ & $\begin{tikzcd}
	\bullet & \bullet & \bullet
	\arrow[from=1-1, to=1-2]
	\arrow[from=1-2, to=1-3]
\end{tikzcd}$

\\

\\

\includegraphics[align=c,scale=1.5]{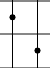} & $\mapsto$ & $\begin{tikzcd}
	\bullet & \bullet & \bullet
	\arrow[""{name=0, anchor=center, inner sep=0}, curve={height=12pt}, from=1-1, to=1-2]
	\arrow[""{name=1, anchor=center, inner sep=0}, curve={height=-12pt}, from=1-1, to=1-2]
	\arrow[""{name=2, anchor=center, inner sep=0}, curve={height=12pt}, from=1-2, to=1-3]
	\arrow[""{name=3, anchor=center, inner sep=0}, curve={height=-12pt}, from=1-2, to=1-3]
	\arrow[shorten <=3pt, shorten >=3pt, Rightarrow, from=1, to=0]
	\arrow[shorten <=3pt, shorten >=3pt, Rightarrow, from=3, to=2]
\end{tikzcd}$

\\

\\

\includegraphics[align=c,scale=1.5]{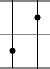} & $\mapsto$ & $\begin{tikzcd}
	\bullet & \bullet & \bullet
	\arrow[""{name=0, anchor=center, inner sep=0}, curve={height=12pt}, from=1-1, to=1-2]
	\arrow[""{name=1, anchor=center, inner sep=0}, curve={height=-12pt}, from=1-1, to=1-2]
	\arrow[""{name=2, anchor=center, inner sep=0}, curve={height=12pt}, from=1-2, to=1-3]
	\arrow[""{name=3, anchor=center, inner sep=0}, curve={height=-12pt}, from=1-2, to=1-3]
	\arrow[shorten <=3pt, shorten >=3pt, Rightarrow, from=1, to=0]
	\arrow[shorten <=3pt, shorten >=3pt, Rightarrow, from=3, to=2]
\end{tikzcd}$

\\

\\

\includegraphics[align=c,scale=1.5]{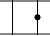} & $\mapsto$ & $\begin{tikzcd}
	\bullet & \bullet & \bullet
	\arrow[""{name=0, anchor=center, inner sep=0}, curve={height=12pt}, from=1-2, to=1-3]
	\arrow[""{name=1, anchor=center, inner sep=0}, curve={height=-12pt}, from=1-2, to=1-3]
	\arrow[from=1-1, to=1-2]
	\arrow[shorten <=3pt, shorten >=3pt, Rightarrow, from=1, to=0]
\end{tikzcd}$

\\

\\

\includegraphics[align=c,scale=1.5]{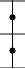} & $\mapsto$ & $\begin{tikzcd}
	\bullet & \bullet
	\arrow[""{name=0, anchor=center, inner sep=0}, from=1-1, to=1-2]
	\arrow[""{name=1, anchor=center, inner sep=0}, curve={height=-18pt}, from=1-1, to=1-2]
	\arrow[""{name=2, anchor=center, inner sep=0}, curve={height=18pt}, from=1-1, to=1-2]
	\arrow[shorten <=2pt, shorten >=2pt, Rightarrow, from=1, to=0]
	\arrow[shorten <=2pt, shorten >=2pt, Rightarrow, from=0, to=2]
\end{tikzcd}$ 
 \end{tabular}

\end{center}

\end{example}

\begin{definition}
	
	Given a simple $k$-string diagram $d$ and a $k$-globular set $X$, we define an \textbf{$X$-labelling} of $d$ to be an $X$-labelling of $\pi(d)$, i.e. a map $\widehat{\pi(d)}\to X$. 
	
\end{definition}

\begin{definition}
 
 We define $T_n^{\sd}:\gSet_n\to\gSet_n$ by $$T_n^{\sd}(X)(k)=\coprod_{d\in\sd(k)}\gSet_n(\widehat{\pi(d)},X).$$ We equip this with the structure of a monad by interpreting a simple string diagram labelled by simple string diagrams as a simple string diagram (see \cite{sesquicat}). An \textbf{$n$-sesquicategory} is a $T_n^{\sd}$-algebra. A \textbf{functor} of $n$-sesquicategories is a morphism of $T_n^{\sd}$-algebras. 
 
\end{definition}

The map $\pi:\sd\to\pd$ now induces a map of monads $$\pi:T_n^{\sd}\to T_n.$$

In \cite{sesquicat} we gave a \textbf{presentation} of the monad $T_n^{\sd}$ by generators and relations. Generators are of the form $\circ_{i,j}$ and $u_i$, for $i,k=1,\cdots,n$. The generators $\circ_{i,j}$ induce binary operations called \textbf{composition} (also called \textbf{whiskering}, if $i\neq j$). The generators $u_i$ induce unary operations which create the \textbf{identity} morphisms. The relations essentially express \textbf{associativity} and \textbf{unitality} of composition. We can therefore say that an $n$-sesquicategory is an $n$-globular set equipped with strictly associative and unital composition operations, which are however not required to satisfy the usual Godement interchange relations which hold in strict $n$-categories.

\begin{remark}

In order to obtain a theory of semistrict $n$-categories, one must enlarge the monad to include cells implementing coherent versions of these interchange relations. This is work in progress.
 
\end{remark}

\begin{example}
 
Here are all generators $\circ_{i,j}$ and $u_i$ for $i,j\leq 4$, as simple string diagrams.

\begin{center}\begin{tabular}{llllll}

$\circ_{1,1}=\includegraphics[align=c]{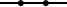}$ & $\circ_{1,2}=\includegraphics[align=c]{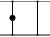}$ & $\circ_{2,1}=\includegraphics[align=c]{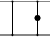}$ & $\circ_{2,2}=\includegraphics[align=c]{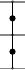}$ & $\circ_{1,3}=\includegraphics[align=c]{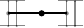}$ &  $\circ_{3,1}=\includegraphics[align=c]{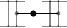}$
 
\end{tabular}\end{center}

\begin{center}\begin{tabular}{lllll}

  $\circ_{2,3}=\includegraphics[align=c]{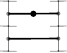}$ & $\circ_{3,2}=\includegraphics[align=c]{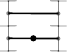}$ & $\circ_{3,3}=\includegraphics[align=c]{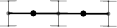}$ & $\circ_{1,4}=\includegraphics[align=c]{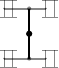}$ & $\circ_{4,1}=\includegraphics[align=c]{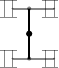}$ 
 
\end{tabular}\end{center}

\begin{center}\begin{tabular}{lllll}

 $\circ_{2,4}=\includegraphics[align=c]{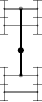}$ & $\circ_{4,2}=\includegraphics[align=c]{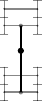}$ & $\circ_{3,4}=\includegraphics[align=c]{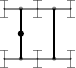}$ &   $\circ_{4,3}=\includegraphics[align=c]{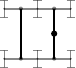}$ & $\circ_{4,4}=\includegraphics[align=c]{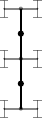}$
 
\end{tabular}\end{center}

\begin{center}\begin{tabular}{llll}

$u_1=\includegraphics[align=c]{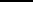}$ & $u_2=\includegraphics[align=c]{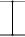}$ & $u_3=\includegraphics[align=c]{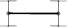}$ & $u_4=\includegraphics[align=c]{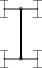}$
 
\end{tabular}\end{center}
 
\end{example}

Finally, in \cite{sesquicat} we also give an \textbf{inductive description} of $n$-sesquicategories similar to the usual inductive definition of strict $n$-category. We won't need it in the present paper.

\subsection{Computads and string diagrams for $n$-sesquicategories}

Both $T_n$ and $T_n^{\sd}$ are finitary monads, so they admit a theory of \textbf{computads}. This is explained in detail in \cite{csp_thesis}[Section 2.7], so we give only a brief summary here. Given a a  monad $T$ on $\gSet_n$ one denotes by $\Alg_T$ its category of algebras. If $T$ is finitary one can inductively define categories $\Comp^T_k$ whose objects are the $k$-computads with respect to $T$, for $0\leq k\leq n+1$, together with adjunctions \[\begin{tikzcd}
	{F_k:\Comp_k^T} & \Alg_T:V_k
	\arrow[""{name=0, anchor=center, inner sep=0}, shift left=2, from=1-1, to=1-2]
	\arrow[""{name=1, anchor=center, inner sep=0}, shift left=2, from=1-2, to=1-1]
	\arrow["\dashv"{anchor=center, rotate=-90}, draw=none, from=0, to=1]
\end{tikzcd}.\] A \textbf{$k$-computad} $C$ consists of sets $C_m$ of generating $m$-cells for each $0\leq m\leq k$, together with source and target maps $s,t:C_m\to F_{m-1}(C)_{m-1}$. We refer to $F_k(C)$ as the $T$-algebra presented by $C$. It is defined from $F_{k-1}(C)$ by a certain pushout in $T$-algebras. 

When $T=T_n$ and $k\leq n$, the morphisms in the strict $n$-category $F_k(C)$ presented by $C$ are the $n$-categorical composites of cells in $C$. We obtain $F_{n+1}(C)$ from $F_n(C)$ by imposing on $n$-morphisms the equivalence relation generated by $s(x)\sim t(x)$ for each $x\in C_{n+1}$. 

\begin{definition}

A $n$-categorical \textbf{presentation} is an $(n+1)$-computad for $T_n$.
 
\end{definition}

\begin{example}\label{computad}
	
	Let's define a $2$-computad $C$ having four $0$-cells, which we call $X, Y, Z$ and $W$, five $1$-cells, $f,g:X\to Z$, $a:X\to Y$, $b:Z\to W$ and $c:Y\to W$ and two $2$-cells, namely $\alpha:f\to g$ and $\beta:b\circ g\to c\circ a$. Then the following composite is a $2$-morphism in $F_2(C)$.
	
\[\begin{tikzcd}
	X & Y \\
	Z & W
	\arrow["a", from=1-1, to=1-2]
	\arrow[""{name=0, anchor=center, inner sep=0}, from=1-1, to=2-1]
	\arrow[""{name=1, anchor=center, inner sep=0}, "f"', curve={height=30pt}, from=1-1, to=2-1]
	\arrow["c", from=1-2, to=2-2]
	\arrow["\beta", Rightarrow, from=2-1, to=1-2]
	\arrow["b"', from=2-1, to=2-2]
	\arrow["\alpha", shorten <=6pt, shorten >=6pt, Rightarrow, from=1, to=0]
\end{tikzcd}: b\circ f \to c\circ a\]
	
\end{example}

The above example shows how computads allow us to make sense of composing pasting diagrams which are more complicated than the globular pasting diagrams appearing in the definition of strict $n$-categories. In \cite{sesqui_comp} we apply the same idea to $n$-sesquicategories and string diagrams. We now give a brief overview.

\begin{notation}

We denote by $\one$ the terminal $(n+1)$-computad for $T_n^{\sd}$.
 
\end{notation}

\begin{definition}

Given an $n$-computad $C$ for $T_n^{D^s}$, a \textbf{$C$-labelled $k$-string diagram} is a $k$-morphism in $F_n(C)$. An \textbf{unlabelled $k$-string diagram} is a $k$-morphism in $F_n(\one)$. The \textbf{shape} of a $C$-labelled string diagram is its image under the map $F_n(C)\to F_n(\one)$. A \textbf{cell shape} is a cell of $\one$.

\end{definition}

\begin{remark}
 
Given an $(n+1)$-computad for $T_n^{\sd}$, the $k$-morphisms in $F_{n+1}(C)$ are the $C$-labelled $k$-string diagrams, for $k\leq n-1$. The $n$-morphisms in $F_{n+1}(C)$ are equivalence classes of $C$-labelled $n$-string diagrams, where the equivalence relation is generated by $(n+1)$-cells in $C$.
 
\end{remark}

Using the presentation of $T_n^{\sd}$ by generators and relations, one can describe $C$-labelled $k$-string diagrams as equivalence classes of trees whose internal vertices are labelled by the generating operations of $T_n^{\sd}$ and whose leaves are labelled by cells in $C$. The equivalence relation on trees is generated by the relations in the presentation of $T_n^{\sd}$ (see \cite{sesqui_comp} for details).

In \cite{sesqui_comp} we introduced a notion of \textbf{normal form} for such labelled trees and we proved that each equivalence class contains a unique tree in normal form. This allows us to construct, for each unlabelled $k$-string diagram $d\in F_n(\one)_k$, a $k$-computad $\hat{d}$ together with a natural bijection between $C$-labelled $k$-string diagrams of shape $d$ and maps $\hat{d}\to C$, for any computad $C$. 

\begin{remark}
 This amounts to saying that the functor $\Comp_{n}^{T_n^{\sd}}\to \Set$ sending an $n$-computad $C$ to the set of $C$-labelled $k$-string diagrams $F_n(C)_k$ is \textbf{familially representable} in the sense of \cite{Carboni_Johnstone_1995}, with representing family $(\hat{d})_{d\in F_n(\one)_k}$.
\end{remark}

\begin{remark}
 
This justifies the name $C$-labelled string diagrams for morphisms in $F_n(C)$. Such a morphism consists exactly of an unlabelled string diagram $d$ (its shape) together with a $C$-labelling (the corresponding map $\hat{d}\to C$). 
 
\end{remark}

\begin{remark}

From this it follows that $\Comp_{n+1}^{T_n^{\sd}}$ is equivalent to the category of presheaves on a small category whose objects are cell shapes and where a morphism $c_1\to c_2$ is given by a map of computads $\widehat{c_1}\to\widehat{c_2}$ (see \cite{sesqui_comp}).

\end{remark}

\begin{remark}

Note that this famously fails for strict $n$-categories (see \cite{3computads_closed},\cite{direct_comp_closed}).
 
\end{remark}

\begin{notation}

Given a $C$-labelling $\lambda:\hat{d}\to C$ of some unlabelled diagram $d$, we denote by $d[\lambda]$ the corresponding $C$-labelled diagram.
 
\end{notation}

Another important consequence of the existence of normal forms is that they make it possible to construct \textbf{graphical representations of unlabelled $k$-string diagrams}. We construct these inductively from the corresponding tree in normal form. We already have a way of representing the generators of $T_n^{\sd}$ which label the internal nodes. The $k$-cells which label the leaves can be represented by drawing the source and target $(k-1)$-diagrams and connecting the dots representing $(k-1)$-cells to a common node in the middle, representing the $k$-cell. Then one composes pictures according to the structure of the tree, going from the leaves to the root (see \cite{sesqui_comp} for more details).

\begin{example}

We give some examples of cell shapes and unlabelled string diagrams. Here is a $2$-cell shape, with its source and target $1$-diagrams.

\begin{center}\begin{tabular}{llclc}
 
$\includegraphics[align=c]{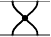}$ & $:$ & $\includegraphics[align=c]{figures/circ11.pdf}$ & $\to$ & $\includegraphics[align=c]{figures/circ11.pdf}$

\end{tabular}\end{center}

Here is a tree in normal  form, then the same tree where we replace each generator by the corresponding picture and finally the picture of the diagram itself.

\[\vcenter{\vbox{\xymatrixcolsep{.6pc}\xymatrixrowsep{1pc}\xymatrix{ & & \circ_{2,2} & & \\ & \circ_{2,1}\ar[ru] & & \circ_{1,2}\ar[lu] & \\ \includegraphics[align=c,scale=.5]{figures/2cell_1.pdf}\ar[ru] & & \includegraphics[align=c,scale=.5]{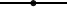}\ar[lu]\text{   }\includegraphics[align=c,scale=.5]{figures/1cell.pdf}\ar[ru] & & \includegraphics[align=c,scale=.5]{figures/2cell_1.pdf}\ar[lu]}}}=\includegraphics[align=c]{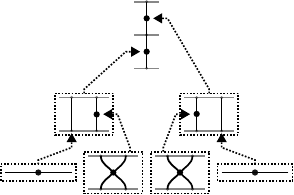}=\includegraphics[align=c]{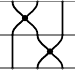}.\]
 Here is a $3$-cell with its source and target $2$-diagrams.

\begin{center}\begin{tabular}{llclc}
 
$\includegraphics[align=c]{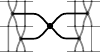}$ & $:$ & $\includegraphics[align=c]{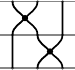}$ & $\to$ & $\includegraphics[align=c]{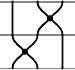}$

\end{tabular}\end{center}

Here is a $4$-cell with its source and target $3$-diagrams.

$$\includegraphics[align=c]{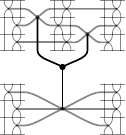}: \includegraphics[align=c]{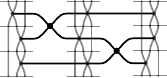}\to\includegraphics[align=c]{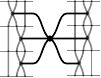}$$
 
\end{example}

The map of monads $\pi:T_n^{D^s}\to T_n$ induces functors \begin{center}$\pi^*:\Alg_{T_n}\to\Alg_{T_n^{\sd}}$ and $\pi_*:\Comp_k^{T_n^{D^s}}\to \Comp_{k}^{T_n}$.\end{center} The computads $C$ and $\pi_*(C)$ have the same sets of generating cells in each dimension. Given a $k$-computad $C$ for $T_n^{D^s}$, we have a map of $T_n^{D^s}$-algebras $F_k(C)\to \pi^*F_k(\pi_*(C)).$

\begin{notation}
 
Given a $k$-computad $C$ for $T_n^{\sd}$, we denote the map above by  $$\pi_C:F_k(C)\to \pi^*F_k(\pi_*(C)).$$ 
 
\end{notation}

\begin{example}
 
Here are the pasting diagrams corresponding to the images under $\pi_{\one}$ of some of the string diagrams above.

\begin{center}
\begin{tabular}{lll}

\includegraphics[align=c,scale=1.5]{figures/2cell_1.pdf} & $\mapsto$ & $\begin{tikzcd}[row sep=small]
	& \bullet \\
	\bullet && \bullet \\
	& \bullet
	\arrow[from=2-1, to=1-2]
	\arrow[from=1-2, to=2-3]
	\arrow[from=2-1, to=3-2]
	\arrow[from=3-2, to=2-3]
	\arrow[shorten <=6pt, shorten >=6pt, Rightarrow, from=1-2, to=3-2]
\end{tikzcd}$

\\

\\

\includegraphics[align=c,scale=1.5]{figures/2diag_1.pdf} & $\mapsto$ &$\begin{tikzcd}[row sep=small]
	& \bullet \\
	\bullet && \bullet \\
	& \bullet && \bullet \\
	&& \bullet
	\arrow[from=2-1, to=1-2]
	\arrow[from=1-2, to=2-3]
	\arrow[from=2-1, to=3-2]
	\arrow[from=3-2, to=2-3]
	\arrow[shorten <=6pt, shorten >=6pt, Rightarrow, from=1-2, to=3-2]
	\arrow[from=2-3, to=3-4]
	\arrow[from=3-2, to=4-3]
	\arrow[shorten <=6pt, shorten >=6pt, Rightarrow, from=2-3, to=4-3]
	\arrow[from=4-3, to=3-4]
\end{tikzcd}$
 
\end{tabular}

\end{center}

\end{example}
 
 \subsection{Equivalences}
In a strict $n$-category, we say that a $k$-morphism $f:x\to y$ is an \textbf{isomorphism} if there exists another $k$-morphism $f:y\to x$ such that $f\circ g=\id_y$ and $g\circ f=\id_x$. We also say that $f$ is \textbf{invertible} and we call $g$ its \textbf{inverse} (one can show that it is unique). However, we are more interested in a weaker version of this, known as \textbf{equivalence}.

\begin{definition}

	Let $\CC$ be a strict $n$-category. An $n$-morphism $f:x\to y$ in $\CC$ is an \textbf{equivalence} if it is an isomorphism. When $k<n$, a $k$-morphism $f:x\to y$ in $\CC$ is an \textbf{equivalence} when there is another $k$-morphism $g:y\to x$ and equivalences $f\circ g\to\id_y$ and $g\circ f\to\id_x$ in $\CC$. We say that $x$ is \textbf{equivalent} to $y$, and write $x\simeq y$, if there is an equivalence $x\to y$. When $f:x\to y$ is an equivalence, we also call it \textbf{weakly invertible} and any morphism $g:y\to x$ such that $f\circ g\simeq \id_y$ and $g\circ f\simeq \id_x$ is called a \textbf{weak inverse} to $f$. When $f$ is a $k$-morphism and an equivalence we also call it a \textbf{$k$-equivalence}.

\end{definition}

\begin{definition}

		An \textbf{$n$-groupoid} is an $n$-category all of whose morphisms are equivalences.

\end{definition}

Finally, we use the following notion of weak equivalence for functors, which coincides with the one in the folk model structure of \cite{folk}.

\begin{definition}

A map of sets $f:X\to Y$ is a \textbf{weak equivalence} if it is a bijection.

Let $n\geq 1$. A functor $F:\CC\to\DD$ between strict $n$-categories is called \textbf{essentially surjective} if for every object $d\in\DD$ there exists an object $c\in\CC$ and an equivalence $F(c)\to d$ in $\DD$. A functor $F:\CC\to\DD$ between strict $n$-categories is called a \textbf{weak equivalence} if it is essentially surjective and for all objects $c_1,c_2\in\CC$ the induced functor $\CC(c_1,c_2)\to\DD(F(c_1),F(c_2))$ is a weak equivalence of $(n-1)$-categories.

\end{definition}
 
\section{String diagrams for strict $4$-categories}\label{stringdiag}
 
Now we introduce the string diagram notation for strict $4$-categories. This will be based on string diagrams for $4$-sesquicategories and the \textbf{homotopy generators} introduced in \cite{data_struct_quasistrict}. Essentially these homotopy generators are certain distinguished generating cells in a computad which implement the interchange laws which hold in strict $4$-categories but not in $4$-sesquicategories.

\subsection{Homotopy generators} 

Now we want to extend our string diagram notation by introducing cells that connect $C$-labelled string diagrams which have the same image under $\pi_C$. In \cite{data_struct_quasistrict} the authors introduce the concept of a \textbf{$4$-signature with homotopy generators} with this same purpose in mind. The notion of a $4$-signature used there corresponds exactly to what we refer to here as a $5$-computad for the monad $T_4^{D^s}$. A $4$-signature with homotopy generators, in the sense of \cite{data_struct_quasistrict}, is then a $5$-computad $C$ for the monad $T_4^{D^s}$ which contains some distinguished cells, called homotopy generators, which the authors describe in detail. The source and target $C$-labelled string diagrams of each homotopy generator have the same image under $\pi_{C}$.

\begin{notation}

In this section, a \textbf{computad} will always mean a $5$-computad (for $T_4^{\sd}$ or $T_4$).

\end{notation}

\begin{remark}
 
The question of whether any two $C$-labelled string diagrams whose images under $\pi_C$ agree can be connected by a sequence of homotopy generators is not addressed in \cite{data_struct_quasistrict}. Proving this would show that the given list of homotopy generators is exhaustive. In an upcoming paper, we address this in the case of $3$-categories. We are also working on extending this to $4$-categories. This is an important step in developing a theory of \textbf{semistrict $4$-categories}, where composition is strictly associative and unital, but the interchange laws hold only up to coherent equivalence. For the purposes of the present paper, however, this issue is not essential. This is because we will be using the string diagram calculus to prove statements about strict $4$-categories, so all we need is a list of homotopy generators which is enough for our computations.

\end{remark}

\begin{remark}
 
Here homotopy generators are somewhat unnaturally added as extra cells to computads. In the theory of semistrict $4$-categories mentioned above, the homotopy generators will come into the definition of a monad $T_4^{ss}$ whose algebras will be semistrict $4$-categories. 
 
\end{remark}

We now present the definition of signature with homotopy generators in our language. We start by describing the \textbf{unlabelled homotopy generators}. These are an infinite collection of unlabelled cells, which come in finitely many \textbf{types}. Each unlabelled homotopy generator $H$ has the property that $\pi_{\one}(s(H))=\pi_{\one}(t(H))$.

We now describe the \textbf{generator types} from \cite{data_struct_quasistrict} in our string diagram notation. For each generator of a certain type, there is another one where the source and target are switched. In the case of $4$-dimensional generators, there are also relations (i.e. $5$-cells) making the pair into an isomorphism. In the case of $3$-dimensional generators there are $4$-cells and relations extending the pair into an adjoint equivalence. For each of the generator types below there are also other versions, which correspond to switching the relative positions of the various components of the source and target. We omit them for conciseness. For more details, see \cite{data_struct_quasistrict}.

The generators of type $\I_2$ are the $3$-cells \[\includegraphics[align=c,scale=1.5]{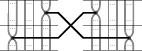}:\includegraphics[align=c,scale=1.5]{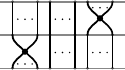}\to \includegraphics[align=c,scale=1.5]{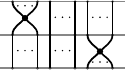}.\] 

The generators of type $\I_3$ are the $4$-cells \[\includegraphics[align=c,scale=1.5]{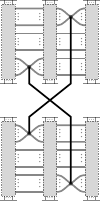}:\includegraphics[align=c,scale=1.5]{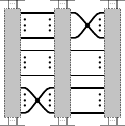}\to \includegraphics[align=c,scale=1.5]{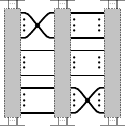}.\] 

The generators of type $\I_4$ are the relations between $4$-diagrams \[\includegraphics[align=c,scale=1.5]{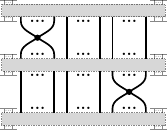}= \includegraphics[align=c,scale=1.5]{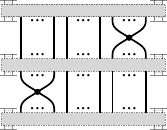} .\] 

We call type $\I_k$ homotopy generators \textbf{interchangers}.

The generators of type $\II_3$ are the $4$-cells  \[\includegraphics[align=c,scale=1.5]{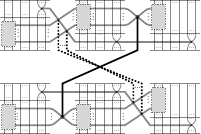}:\includegraphics[align=c,scale=1]{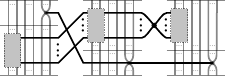}\to \includegraphics[align=c,scale=1]{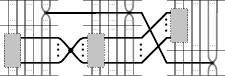}\]

The generators of type $\II_4$ are the relations between $4$-diagrams  \[\includegraphics[align=c,scale=1]{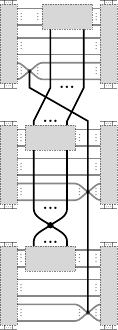}= \includegraphics[align=c,scale=1]{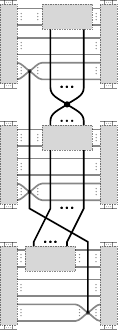}.\] 

We call type $\II_k$ homotopy generators \textbf{pull-through} cells.

The homotopy generators of type $\III_4$ are the relations between $4$-diagrams \[\includegraphics[align=c,scale=1]{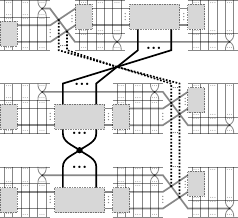}= \includegraphics[align=c,scale=1]{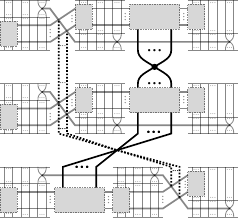}.\]

The homotopy generators of type $\IV_4$ are the relations between $4$-diagrams  \[\includegraphics[align=c,scale=1]{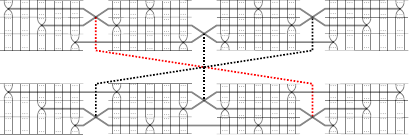}= \includegraphics[align=c,scale=1]{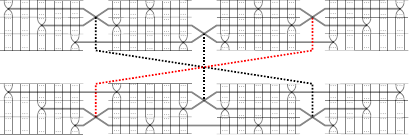}.\] Here the $3$-diagrams consist only of type $\I_2$ homotopy generators. The $4$-diagrams are two different type $\II_3$ homotopy generators, each applying pull through to an interchanger (highlighted in red).

The homotopy generators of type $\V_4$ are the relations between $4$-diagrams  \[\includegraphics[align=c,scale=1]{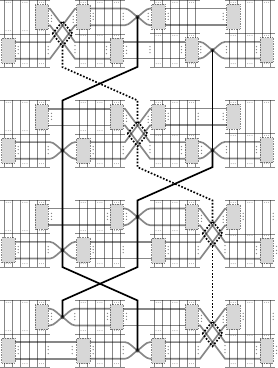}= \includegraphics[align=c,scale=1]{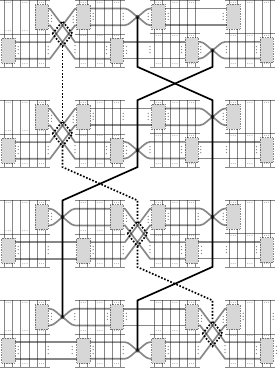}.\]

The homotopy generators of type $\VI$ are the relations between $4$-diagrams  \[\includegraphics[align=c,scale=1]{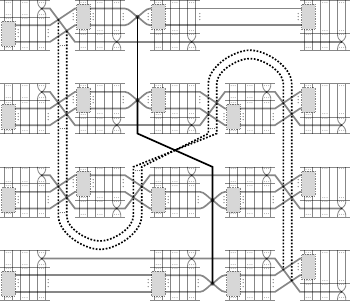}= \includegraphics[align=c,scale=1]{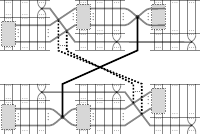}.\]

\begin{remark}

For each homotopy generator $H$, we have $\pi_{\one}(s(H))=\pi_{\one}(t(H))$ and also $\widehat{s(H)}=\widehat{t(H)}$.
 
\end{remark}

\begin{definition}
 
A \textbf{signature with homotopy generators} is a computad $C$ for $T_4^{\sd}$ such that, for each unlabelled homotopy generator $H$ of dimension $k$ and each $C$-labelling $\lambda:\widehat{s(H)}=\widehat{t(H)}\to C$ there is a distingished $(k+1)$-cell $H[\lambda]:s(H)[\lambda]\to t(H)[\lambda]$ in $C$. 
 
\end{definition}

\subsection{String diagrams in strict $4$-categories}

We now explain how we will use string diagrams over signatures with homotopy generators to do calculations in strict $4$-categories. 

Given a signature $\Sigma\in\Comp_5^{T_4^{\sd}}$ with homotopy generators, define $\bar{\pi}_*(\Sigma)\in\Comp_5^{T_4}$ to be the computad for $T_4$ obtained from $\pi_*(\Sigma)$ by discarding the cells corresponding to the homotopy generators and replacing these cells by identities, whenever they appear in a source or target diagram for another cell.

Now let $\CC$ be a strict $4$-category and consider a functor $F_5(\bar{\pi}_*(\Sigma))\to \CC$. Such a functor consists of a choice of $k$-morphism in $\CC$ for each $k$-cell in $\Sigma$ which is not a homotopy generator, such that the source and target of each chosen morphism is a prescribed composite of chosen morphisms in lower dimension and the $4$-morphisms satisfy relations corresponding to the $5$-cells in $\Sigma$. We can extend this map to a map $F_5(\pi_*(\Sigma))\to \CC$ by sending all homotopy generators to identity morphisms in $\CC$. 

We pull back to a map $\pi^*F_5(\pi_*(\Sigma))\to \pi^*\CC$  and compose with $F_5(\Sigma)\to \pi^*F_5(\pi_*(\Sigma))$ to get a map of $T_4^{D^s}$-algebras $F_5(\Sigma)\to \pi^*\CC$, so any string diagram in the signature $\Sigma$ determines a morphism in $\CC$ and any identity between $4d$ string diagrams in the signature $\Sigma$ implies an identity between the corresponding $4$-morphisms in $\CC$. This is what allows us to use string diagrams with labels in a strict $4$-category to describe composite morphisms in this $4$-category. 

\begin{remark}
 
The string diagram calculus we are describing would be most naturally interpreted in some kind of semistrict $4$-category. The procedure described here collapses all the homotopies and so interprets this in a strict $4$-category. 
 
\end{remark}

\section{Functor $4$-categories}\label{functorcat}

Given $n$-categories $\CC$ and $\DD$, one can define an $n$-category $\Fun(\CC,\DD)$, whose $k$-morphisms are called \textbf{$k$-transfors}. A $0$-transfor is a functor, a $1$-transfor is a natural transformation, a $2$-transfor is also called a modification and a $3$-transfor is sometimes known as a perturbation (see the nLab page "transfor" for a discussion of this terminology).

Using the left and right internal $\Hom$ from the monoidal biclosed structure on $n$-categories associated to the Crans-Gray tensor product (\cite{crans95}) one can define $n$-categories $\Fun^{\lax}(\CC,\DD)$ and $\Fun^{\oplax}(\CC,\DD)$ for strict $n$-categories $\CC$ and $\DD$. One can check that a $k$-morphism in $\Fun^{\oplax}(\CC,\DD)$ is a rule that associates to each $\ell$-morphism in $\CC$ a map $\theta^{(k);(\ell)}\to\DD$, satisfying certain relations of compatibility with composition. Here $\theta^{(k);(\ell)}$ is the $(k+\ell)$-computad explictly constructed in \cite{freyd_scheim}. It can also be described as the Crans-Gray tensor product $\theta^{(k)}\otimes\theta^{(\ell)}$, where $\theta^{(k)}$ denotes the computad generated by a single $k$-cell. Similarly, a $k$-morphism in $\Fun^{\lax}(\CC,\DD)$ is a rule that associates to each $\ell$-morphism in $\CC$ a map $\theta^{(\ell);(k)}\to\DD$.

One can then define the $n$-category $\Fun(\CC,\DD)$ as the subcategory of $\Fun^{\oplax}(\CC,\DD)$ consisting of those $k$-morphisms which associate to an $\ell$-morphism in $\CC$ a $(k+\ell)$-equivalence in $\DD$, for $k,\ell\geq 1$. Below we give an explicit description of $\Fun(\CC,\DD)$ in terms of string diagrams, when $\CC$ and $\DD$ are $4$-categories. 

\begin{remark}
 
Whereas in the previous section we were using dashed lines to denote interchangers in lower dimension, here these are omitted from the notation, for simplicity. We will instead use dashed lines as well as colours to denote labellings of the diagrams.
 
\end{remark}

\subsection{Natural transformations}

Given functors $F,G:\CC\to\DD$, a \textbf{natural transformation}, or $1$-transfor, $\alpha:F\to G$ consists of the following data. We use {\red} and {\blue} to denote the images of objects and morphisms under $F$
and $G$, respectively.
\begin{enumerate}
 \item[0.] For each object $Y\in \CC$ a $1$-morphism $\alpha_Y:F(Y)\to G(Y)$. \begin{center}\begin{tabular}{lcr}$Y=$ \includegraphics[scale=1.5,align=c]{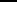} & $\mapsto$ & $\alpha_Y=$  \includegraphics[scale=1.5,align=c]{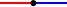}\end{tabular}\end{center}
 
 \item For each $1$-morphism $g:X\to Y$ in $\CC$ an invertible $2$-morphism $\alpha_g$ in $\DD$. \begin{center}\begin{tabular}{lcr} $g=$ \includegraphics[scale=1.5,align=c]{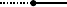} & $\mapsto$ & $\alpha_g=$ \includegraphics[scale=1.5,align=c]{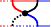} $:$ \includegraphics[scale=1.5,align=c]{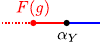} $\to$ \includegraphics[scale=1.5,align=c]{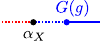}\end{tabular}\end{center}
 
 \item For each $2$-morphism $\zeta:f\to g$ in $\CC$ an invertible $3$-morphism $\alpha_{\zeta}$ in $\DD$. \begin{center}\begin{tabular}{lcr} $\zeta=$ \includegraphics[scale=1.5,align=c]{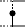} & $\mapsto$ & $\alpha_{\zeta}=$ \includegraphics[scale=1.5,align=c]{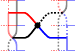} $:$ \includegraphics[scale=1.5,align=c]{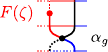} $\to$ \includegraphics[scale=1.5,align=c]{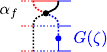}\end{tabular} \end{center}
 
 \item For each $3$-morphism $t:\eta\to\zeta$ in $\CC$ an invertible $4$-morphism $\alpha_{t}$ in $\DD$. \begin{center}\begin{tabular}{lcr}$t=$ \includegraphics[scale=1.5,align=c]{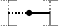} & $\mapsto$ & $\alpha_{t}=$ \includegraphics[scale=1.5,align=c]{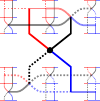} $:$ \includegraphics[scale=1.5,align=c]{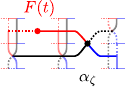} $\to$ \includegraphics[scale=1.5,align=c]{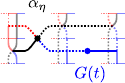}\end{tabular}\end{center}
 
 \item For each $4$-morphism $\W:s\to t$ in $\CC$ a relation $\alpha_{\W}$ in $\DD$.  \begin{center}\begin{tabular}{lcr} $\W=$ \includegraphics[scale=1.5,align=c]{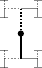} & $\mapsto$ & $\alpha_{\W}:$\includegraphics[scale=1.5,align=c]{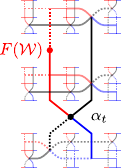} $=$ \includegraphics[scale=1.5,align=c]{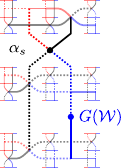}\end{tabular}\end{center}
 
\end{enumerate}

This data is subject to relations equating the values of $\alpha$ on composite morphisms with the corresponding composites of values of $\alpha$ given by stacking diagrams.

\subsection{Modifications}
 
Given natural transformations $\alpha,\beta:F\to G$, a \textbf{modification}, or $2$-transfor, $m:\alpha\to\beta$ consists of the following data. We use {\green} for $\alpha$ and {\purple} for $\beta$. 

\begin{enumerate}
 \item[0.] For each object $Y\in\CC$ a $2$-morphism $m_Y:\alpha_Y\to \beta_Y$ in $\DD$. \begin{center}\begin{tabular}{lcr} $Y=$ \includegraphics[scale=1.5,align=c]{functorcat/1morph/y.pdf} & $\mapsto$ & $m_Y=$ \includegraphics[scale=1.5,align=c]{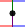} $:$ \includegraphics[scale=1.5,align=c]{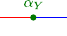} $\to$ \includegraphics[scale=1.5,align=c]{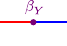}\end{tabular}\end{center}
 
 \item For each $1$-morphism $g:X\to Y$ in $\CC$ an invertible $3$-morphism $m_g$ in $\DD$. \begin{center}\begin{tabular}{lcr}$g=$ \includegraphics[scale=1.5,align=c]{morphisms/f.pdf} & $\mapsto$ & $m_g=$ \includegraphics[scale=1.5,align=c]{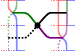} $:$ \includegraphics[scale=1.5,align=c]{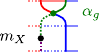} $\to$ \includegraphics[scale=1.5,align=c]{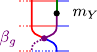}\end{tabular}\end{center}
 
 \item For each $2$-morphism $\zeta=\includegraphics[scale=1.5,align=c]{morphisms/eta.pdf}:f\to g$ in $\CC$ an invertible $4$-morphism $m_{\zeta}$ in $\DD$.  $$m_{\zeta}= \includegraphics[scale=1.5,align=c]{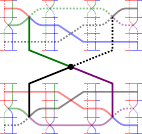}: \includegraphics[scale=1.5,align=c]{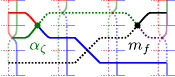}\to \includegraphics[scale=1.5,align=c]{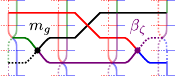}$$
 
 \item For each $3$-morphism $t:\eta\to\zeta$ in $\CC$ a relation $m_t$ in $\DD$. \begin{center}\begin{tabular}{lcr}$t=$ \includegraphics[scale=1.5,align=c]{morphisms/s.pdf} & $\mapsto$ & $m_{t}:$ \includegraphics[scale=1.5,align=c]{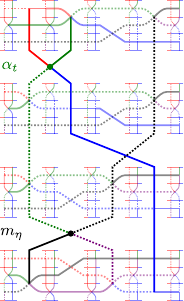} $=$ \includegraphics[scale=1.5,align=c]{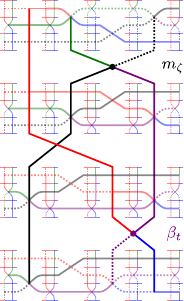}\end{tabular}\end{center}
 
\end{enumerate}

This data is subject to relations equating the values of $m$ on composite morphisms with the corresponding composites of values of $m$ given by stacking diagrams.

\subsection{Perturbations}
 
Given modifications $l,m:\alpha\to\beta$, a \textbf{perturbation}, or $3$-transfor, $\A:l\to m$ consists of the following data. We use {\orange} for $l$ and {\lightblue} for $m$. 

\begin{enumerate}
 \item[0.] For each object $Y\in\CC$ a $3$-morphism $\A_Y:l_Y\to m_Y$ in $\DD$. \begin{center}\begin{tabular}{lcr} $Y=$ \includegraphics[scale=1.5,align=c]{functorcat/1morph/y.pdf} & $\mapsto$ & $\A_Y=$\includegraphics[scale=1.5,align=c]{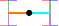} $:$ \includegraphics[scale=1.5,align=c]{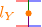} $\to$ \includegraphics[scale=1.5,align=c]{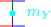}\end{tabular}\end{center}
 
 \item For each $1$-morphism $g=\includegraphics[scale=1.2,align=c]{morphisms/f.pdf}:X\to Y$ in $\CC$ an invertible $4$-morphism $\A_g$ in $\DD$. \begin{center} $\A_g=$ \includegraphics[scale=1.5,align=c]{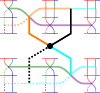} $:$ \includegraphics[scale=1.5,align=c]{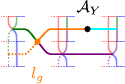} $\to$ \includegraphics[scale=1.5,align=c]{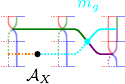}\end{center}
 
 \item For each $2$-morphism $\zeta:f\to g$ in $\CC$ a relation $\A_{\zeta}$ in $\DD$. \begin{center}\begin{tabular}{lcr} $\zeta=$ \includegraphics[scale=1.5,align=c]{morphisms/eta.pdf} & $\mapsto$ & $\A_{\zeta}:$ \includegraphics[scale=1.5,align=c]{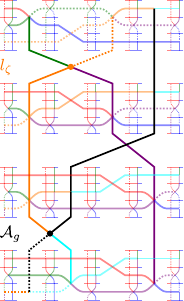} $=$ \includegraphics[scale=1.5,align=c]{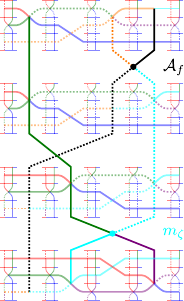}\end{tabular}\end{center}
 
\end{enumerate}

This data is subject to relations equating the values of $\A$ on composite morphisms with the corresponding composites of values of $\A$ given by stacking diagrams. 

\subsection{$4$-transfors}

Given perturbations $\A,\B:l\to m$, a \textbf{$4$-transfor} $\Lambda:\A\to\B$ consists of the following data. We use {\olive} for $\A$ and {\teal} for $\B$.  

\begin{enumerate}
 \item[0.] For each object $Y\in\CC$ a $4$-morphism $\Lambda_Y:\A_Y\to \B_Y$ in $\DD$.  \begin{center}\begin{tabular}{lcr} $Y=$ \includegraphics[scale=1.5,align=c]{functorcat/1morph/y.pdf} & $\mapsto$ & $\Lambda_Y=$ \includegraphics[scale=1.5,align=c]{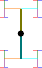} $:$ \includegraphics[scale=1.5,align=c]{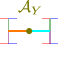} $\to$ \includegraphics[scale=1.5,align=c]{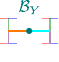}\end{tabular}\end{center}
 
 \item For each $1$-morphism $g=:X\to Y$ in $\CC$ a relation $\Lambda_g$ in $\DD$. \begin{center}\begin{tabular}{lcr}$g=$  \includegraphics[scale=1.5,align=c]{morphisms/f.pdf} & $\mapsto$ & $\Lambda_g:$ \includegraphics[scale=1.5,align=c]{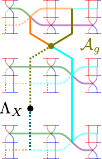} $=$ \includegraphics[scale=1.5,align=c]{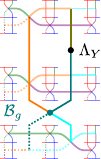}\end{tabular}\end{center}
 
\end{enumerate}

\begin{definition}
Let $\CC, \DD$ be strict $4$-categories. The \textbf{functor category} $\Fun(\CC,\DD)$ is the strict $4$-category whose set of $k$-morphisms is the set of $k$-transfors. Source and target maps are already implicit in the above definitions. It is also clear how to define composition operations by stacking diagrams, which provides the $T_4$-monad structure on $\Fun(\CC,\DD)$. 
\end{definition}

\begin{notation}
	
	We use $\Map$ to denote the underlying $4$-groupoid in $\Fun$, whose $k$-morphisms are weakly invertible $k$-transfors.
	
\end{notation}

\begin{notation}
 
If $\PP$ is a presentation, we write $\Map(\PP,\CC)$ instead of $\Map(F(\PP),\CC)$.
 
\end{notation}

\section{Composition with weakly invertible morphisms}\label{scompfun}

The goal of this section is to give a proof of the following Proposition, which will be the main ingredient in the proof of Theorem \ref{main} on fibrations of mapping $4$-groupoids.

\begin{proposition}\label{compfun}
	
Take $1 \leq k \leq n\leq 4$  and let $\CC$ be an $n$-category. Then any functor between $(n-k)$-categories of $k$-morphisms in $\CC$ given by composition with weakly invertible morphisms of dimension $\leq k$ is essentially surjective. When $k=n$, any such functor is a bijection. 
	
\end{proposition}	

 For the rest of the section, let $\CC$ be a $4$-category, $A,B,C$ objects in $\CC$ and $f:B\to C$ a weakly invertible $1$-morphism in $\CC$. We write $A=$ \includegraphics[scale=1.5,align=c]{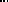}, $B=$ \includegraphics[scale=1.5,align=c]{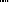}, $C=$ \includegraphics[scale=1.5,align=c]{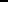} and $f=$ \includegraphics[scale=1.5,align=c]{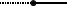}. 
 
 \subsection{Higher invertibility data}
 
 Since $f$ is weakly invertible, we can pick a weak inverse $f^{-1}$, together with equivalences $f^{-1}\circ f \simeq \id_B$ and $f\circ f^{-1} \simeq \id_C$. Now the weakly invertible $2$-morphisms witnessing these two equivalences themselves have weak inverses and there are weakly invertible $3$-morphisms between composites of these $2$-morphisms and the appropriate identity $2$-morphisms and so on. In this way, simply by unraveling the definition of a weak inverse $1$-morphism in a $4$-category, we obtain a collection of $k$-morphisms $(1\leq k\leq 4)$ and some relations between the $4$-morphisms, which we call \textbf{invertibility data} for $f$. We write down here only those higher morphisms and relations that we will use explicitly in the rest of the section. 

\textbf{$1$-morphisms}: 

\begin{longtable}{llcccrr}

	$f=$ \includegraphics[scale=1.5,align=c]{compfun/f.pdf} & $:$ & \includegraphics[scale=1.5,align=c]{compfun/b.pdf} & $\xymatrix@1{\ar@1@<1ex>[r] & \ar@1@<1ex>[l]}$ & \includegraphics[scale=1.5,align=c]{compfun/c.pdf} & $:$ & \includegraphics[scale=1.5,align=c]{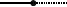} $=f^{-1}$

\end{longtable}

\textbf{$2$-morphisms}: 

\begin{longtable}{rlcccrl}
	
	$u^{-1}=$ \includegraphics[scale=1.5,align=c]{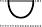} & $:$ & \includegraphics[scale=1.5,align=c]{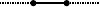} & $\xymatrix@1{\ar@1@<1ex>[r] & \ar@1@<1ex>[l]}$ & \includegraphics[scale=1.5,align=c]{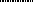} & $:$ & \includegraphics[scale=1.5,align=c]{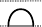} $=u$
	
\\

\\
		
		$c=$ \includegraphics[scale=1.5,align=c]{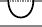} & $:$ & \includegraphics[scale=1.5,align=c]{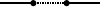} & $\xymatrix@1{\ar@1@<1ex>[r] & \ar@1@<1ex>[l]}$ & \includegraphics[scale=1.5,align=c]{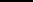} & $:$ & \includegraphics[scale=1.5,align=c]{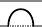} $=c^{-1}$
		
\end{longtable}

\textbf{$3$-morphisms}: 

\begin{longtable}{llcccrr}

\includegraphics[scale=1.5,align=c]{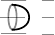} & $:$ & \includegraphics[scale=1.5,align=c]{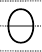} & $\xymatrix@1{\ar@1@<1ex>[r] & \ar@1@<1ex>[l]}$ & 	\includegraphics[scale=1.5,align=c]{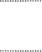} & $:$ & \includegraphics[scale=1.5,align=c]{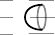}
			
\\

\\
		
\includegraphics[scale=1.5,align=c]{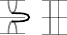} & $:$ & \includegraphics[scale=1.5,align=c]{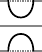} & $\xymatrix@1{\ar@1@<1ex>[r] & \ar@1@<1ex>[l]}$ & 	\includegraphics[scale=1.5,align=c]{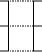} & $:$ & \includegraphics[scale=1.5,align=c]{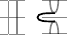}
		
\end{longtable}

\textbf{$4$-morphisms}: 

\begin{longtable}{llcccrr}

\includegraphics[scale=1.5,align=c]{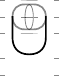} & $:$ & \includegraphics[scale=1.5,align=c]{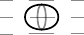} & $\xymatrix@1{\ar@1@<1ex>[r] & \ar@1@<1ex>[l]}$ & 	\includegraphics[scale=1.5,align=c]{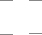} & $:$ & \includegraphics[scale=1.5,align=c]{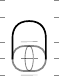}
		
\\

\\
		
\includegraphics[scale=1.5,align=c]{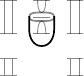} & $:$ & \includegraphics[scale=1.5,align=c]{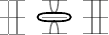} & $\xymatrix@1{\ar@1@<1ex>[r] & \ar@1@<1ex>[l]}$ & 	\includegraphics[scale=1.5,align=c]{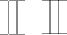} & $:$ & \includegraphics[scale=1.5,align=c]{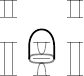}
		
\end{longtable}

\textbf{relations}:

\begin{longtable}{lcrclcr}
		
		\includegraphics[scale=1.5,align=c]{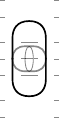} & $=$ & 	\includegraphics[scale=1.5,align=c]{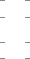} & and &	\includegraphics[scale=1.5,align=c]{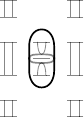} & $=$ & 	\includegraphics[scale=1.5,align=c]{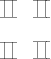} 
		
	\end{longtable}  
	
\subsection{Adjoint equivalences}

We will need to use certain higher morphisms and relations in addition to the ones mentioned above. These don't come simply from the definition of an equivalence in a $4$-category, but from the fact that such an equivalence can be promoted to an adjoint equivalence in $h_2(\CC)$.

\begin{definition}

An \textbf{adjunction} in a $2$-category $\CC$ consists of $1$-morphisms $l:B\to C$ and $r:C\to B$ together with $2$-morphisms $u:\id_B\to r\circ l$ and $c:l\circ r \to \id_C$ satisfying the snake relations. We write $l\dashv r$.

\end{definition} 	

\begin{lemma}
	
If $f:B\to C$ and $f^{-1}:C\to B$ are weakly inverse $1$-morphisms in a $2$-category $\CC$, then there exist invertible $2$-morphisms $u:\id_B\to f^{-1}\circ f$ and $c:f\circ f^{-1} \to \id_C$ satisfying the snake relations, so that we have an adjunction $f\dashv f^{-1}$.	
	
\end{lemma}

\begin{proof}
	
This is well known. There is a nice string diagram proof in the nLab article	on adjoint equivalence.
	
\end{proof}

Now returning to the context of the present section, the pair of weakly inverse $1$-morphisms $(f,f^{-1})$ in the $4$-category $\CC$ can be promoted to an adjoint equivalence in the homotopy $2$-category $h_2(\CC)$, which means that we can pick the $2$-morphisms $u$ and $c$ in the higher invertibility data in such a way that they satisfy snake relations up to equivalence, witnessing the adunction $f\dashv f^{-1}$. This allows us to find pairs of weakly inverse $3$-morphisms witnessing the snake relations, as well as pairs of inverse $4$-morphisms witnessing the fact that these $3$-morphisms are in fact weakly inverse to each other. Of the resulting data, we write down here only the part that we will use explicitly in the rest of the section. Note that this data is not sufficient to specify a fully coherent adjoint equivalence, or a biadjoint biequivalence as it is referred to in \cite{gurski}, in the case of a tricategory. 
	
\textbf{$3$-morphisms}: 	
 
\begin{longtable}{llcccrr}
		
\includegraphics[scale=1.5,align=c]{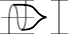} & $:$ & \includegraphics[scale=1.5,align=c]{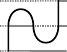} & $\xymatrix@1{\ar@1@<1ex>[r] & \ar@1@<1ex>[l]}$ & 	\includegraphics[scale=1.5,align=c]{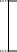} & $:$ & \includegraphics[scale=1.5,align=c]{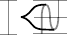}
		
\end{longtable} 

\textbf{$4$-morphisms}: 	
 
\begin{longtable}{llcccrr}
		
\includegraphics[scale=1.5,align=c]{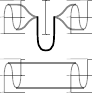} & $:$ & \includegraphics[scale=1.5,align=c]{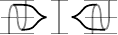} & $\xymatrix@1{\ar@1@<1ex>[r] & \ar@1@<1ex>[l]}$ & 	\includegraphics[scale=1.5,align=c]{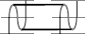} & $:$ & \includegraphics[scale=1.5,align=c]{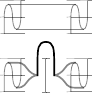}
		
\end{longtable} 

\textbf{relation}:

\begin{longtable}{lcr}
		
		\includegraphics[scale=1.5,align=c]{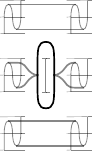} & $=$ & 	\includegraphics[scale=1.5,align=c]{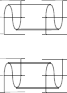}
		
\end{longtable} 

\subsection{Proof of Proposition \ref{compfun}}
	
We now prove a series of Lemmas which taken together imply Proposition \ref{compfun}. We will use a slight extension of our string diagram notation to allow us to specify functors between categories of morphisms. Namely, when we draw a diagram containing a cell labelled by a blank square, this is meant to represent a functor that takes as input a morphism of the appropriate dimension and having the correct source and target, and outputs the result of composig the diagram with this morphism in place of the blank square. 

The proofs of the Lemmas in the rest of this section all follow the same line of reasoning, differing only in the explicit string diagram calculations. We have included all of them for completeness and also because they serve as the first illustrations of the method of string diagram calculations in this paper, but the reader can safely skip most of them.

 Consider the functor $f_*:\Hom(A,B)\to\Hom(A,C)$ given by composition with $f$, which we write down as follows, using the notation we just mentioned: \begin{center}$f_*=$ \includegraphics[scale=1.5,align=c]{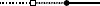} $:\Hom(\includegraphics[scale=1.5,align=c]{compfun/a.pdf},\includegraphics[scale=1.5,align=c]{compfun/b.pdf})\to\Hom(\includegraphics[scale=1.5,align=c]{compfun/a.pdf},\includegraphics[scale=1.5,align=c]{compfun/c.pdf})$.\end{center} 

\begin{lemma}
	
The functor $f_*:\Hom(A,B)\to\Hom(A,C)$ given by composition with $f$	is essentially surjective.
	
\end{lemma}

\begin{proof}
	
Given $g:A\to C$ we have $g\simeq f_*(f^{-1}\circ g)$ for any choice of inverse $f^{-1}$.
	
\end{proof}

Now take $1$-morphisms $a,b:A\to B$ in $\CC$ and consider the functor $f_*:2\Hom(a,b)\to 2\Hom(f\circ a,f\circ b)$ given by composition with $f$. Using {\red} for $a$, {\blue} for $b$ we denote this functor by \begin{center}$f_*=$ \includegraphics[scale=2,align=c]{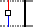} $:2\Hom(\includegraphics[scale=1.5,align=c]{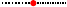},\includegraphics[scale=1.5,align=c]{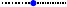})\to 2\Hom(\includegraphics[scale=1.5,align=c]{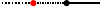},\includegraphics[scale=1.5,align=c]{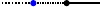})$.\end{center}

\begin{lemma}
	
The functor $f_*:2\Hom(a,b)\to 2\Hom(f\circ a,f\circ b)$ given by composition with $f$ is essentially surjective.
	
\end{lemma}

\begin{proof}

Given a $2$-morphism $$\includegraphics[scale=2,align=c]{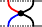}:\includegraphics[scale=2,align=c]{compfun/af.pdf}\to\includegraphics[scale=2,align=c]{compfun/bf.pdf}$$ we need to show that it is equivalent to something in the image of $f_*$. Consider the $2$-morphism $$\includegraphics[scale=2,align=c]{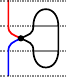}:\includegraphics[scale=2,align=c]{compfun/aa.pdf}\to\includegraphics[scale=2,align=c]{compfun/bb.pdf}.$$ Applying $f_*$ we get \begin{center} \begin{tabular}{cccccc}\includegraphics[scale=2,align=c]{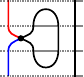} & $\simeq$ & \includegraphics[scale=2,align=c]{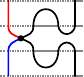} & $\simeq$ & \includegraphics[scale=2,align=c]{compfun/2morph.pdf} & .\end{tabular}\end{center}

\end{proof}

Now take $2$-morphisms $\alpha,\beta:a\to b$ in $\CC$ and consider the functor $f_*:3\Hom(\alpha,\beta)\to 3\Hom(f\circ \alpha,f\circ \beta)$ given by composition with $f$ along $B$. If we use {\green} for $\alpha$ and {\purple} for $\beta$, we can denote $f_*$ by \begin{center}$f_*=$ \includegraphics[scale=2,align=c]{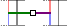} $:3\Hom\left(\text{ }\includegraphics[scale=2,align=c]{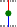},\includegraphics[scale=2,align=c]{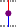}\text{ }\right)\to 3\Hom\left(\text{ } \includegraphics[scale=2,align=c]{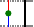},\includegraphics[scale=2,align=c]{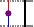} \text{ }\right)$.\end{center} 

\begin{lemma}
	
The functor $f_*:3\Hom(\alpha,\beta)\to 3\Hom(f\circ \alpha,f\circ \beta)$ given by composition with $f$ is essentially surjective.
	
\end{lemma}

\begin{proof}

Given a $3$-morphism $$\includegraphics[scale=2,align=c]{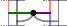}:\includegraphics[scale=2,align=c]{compfun/aaf.pdf}\to\includegraphics[scale=2,align=c]{compfun/bbf.pdf}$$ we need to show that it is isomorphic to something in the image of $f_*$.

 Consider the $3$-morphism $$\includegraphics[scale=2,align=c]{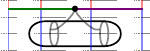}:\includegraphics[scale=2,align=c]{compfun/aaa.pdf}\to\includegraphics[scale=2,align=c]{compfun/bbb.pdf}.$$ Applying $f_*$ we get 

   \begin{longtable}{ll} 
 
\includegraphics[scale=2]{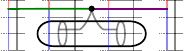} & $\simeq$

\\

\\ 
 
\includegraphics[scale=2]{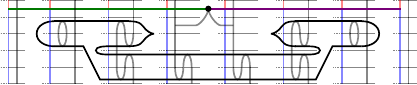} & $\simeq$ 

\\

\\

\includegraphics[scale=2]{compfun/3morph.pdf}. & 

\end{longtable} 

\end{proof}

Finally, take $3$-morphisms $\A,\B:\alpha\to\beta$ and consider the functor $f_*:4\Hom(\A,\B)\to 4\Hom(f\circ \A,f\circ \B)$ given by composition with $f$. If we use {\orange} for $\A$ and {\lightblue} for $\B$, we can denote $f_*$ by  $$f_*=\includegraphics[scale=2,align=c]{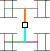}:4\Hom\left(\text{ }\includegraphics[scale=2,align=c]{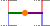},\includegraphics[scale=2,align=c]{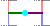}\right)\to 4\Hom\left( \includegraphics[scale=2,align=c]{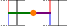},\includegraphics[scale=2,align=c]{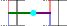}\right).$$ 

\begin{lemma}

The functor $f_*:4\Hom(\A,\B)\to 4\Hom(f\circ \A,f\circ \B)$ given by composition with $f$ is essentially surjective.	
	
\end{lemma}

\begin{proof}

 Given a $4$-morphism $f\circ \A\to f\circ \B$, which we denote by $$\includegraphics[scale=2,align=c]{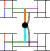}:\includegraphics[scale=2,align=c]{compfun/aaaf.pdf}\to\includegraphics[scale=2,align=c]{compfun/bbbf.pdf},$$ we need to show that it is in the image of $f_*$. 

 Consider the $4$-morphism $$\includegraphics[scale=2,align=c]{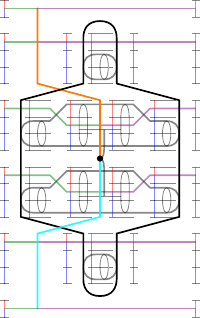}:\includegraphics[scale=2,align=c]{compfun/aaaa.pdf}\to\includegraphics[scale=2,align=c]{compfun/bbbb.pdf}.$$ Applying $f_*$ we get
 
\begin{longtable}{llll}

\includegraphics[scale=1.3]{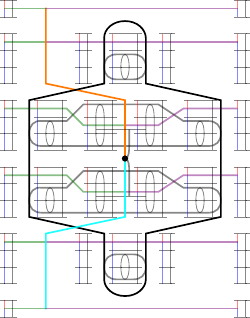} & $=$ & \includegraphics[scale=1.3]{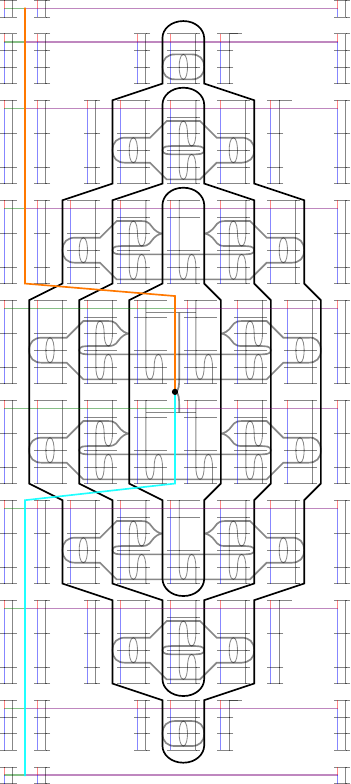} & $=$

\end{longtable}

\begin{longtable}{llll}

$=$ & \includegraphics[scale=1.5]{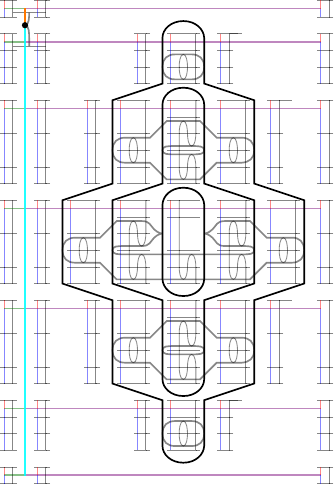} & $=$ & \includegraphics[scale=1.5]{compfun/4morph.pdf}.  

\end{longtable}

\end{proof}

These four Lemmas taken together give a proof of the following Lemma.

\begin{lemma}

Take $1 \leq k \leq n\leq 4$  and let $\CC$ be an $n$-category. Then any functor between $(n-k)$-categories of $k$-morphisms in $\CC$ given by composition with a fixed weakly invertible $1$-morphism is essentially surjective.

\end{lemma}	

By replacing $\CC$ with the appropriate categories of morphisms in $\CC$, we get the following Lemma.

\begin{lemma}
	
Take $1 \leq \ell \leq k \leq n\leq 4$  and let $\CC$ be an $n$-category. Then any functor between $(n-k)$-categories of $k$-morphisms in $\CC$ given by composition with a fixed weakly invertible $\ell$-morphism along an $(\ell-1)$-morphism is essentially surjective.
	
\end{lemma}	
	
In the proofs of these Lemmas, we have actually constructed maps $$k\Hom(f\circ x,f\circ y)\to k\Hom(x, y),$$ where $x$ and $y$ are some $k$-morphisms. When $k=n$, we have showed that this map is right inverse to $f_*$. But it is easy to see that it is also left inverse to $f_*$, which proves the following Lemma.

\begin{lemma}
	
Take $1 \leq \ell \leq n\leq 4$  and let $\CC$ be an $n$-category. Then any functor between sets of $n$-morphisms in $\CC$ given by composition with a fixed weakly invertible $\ell$-morphism along an $(\ell-1)$-morphism is a bijection.
	
\end{lemma}	

\begin{proof}[of Proposition \ref{compfun}]

Any such functor can be obtained as a composite of functors, each of which is given by composition with a fixed weakly invertible $\ell$-morphism ($\ell\leq k$) along an $(\ell-1)$-morphism. A composite of essentially surjective functors is essentially surjective, and a composite of bijective maps is bijective, so we have proved the Proposition.

\end{proof}

\section{Fibrations of mapping $4$-groupoids}\label{sfibration}

The goal of this section is to give a proof of the following Theorem.

\begin{theorem}\label{fibration}

Let $\CC$ a strict $4$-category and $\PP, \QQ$ presentations, with $\QQ$ obtained from $\PP$ by adding a finite number of cells. Then the restriction map $$\Map(\QQ,\CC)\to\Map(\PP,\CC)$$ is a fibration of $4$-groupoids.

\end{theorem}

\begin{definition}
	
	A map of $n$-groupoids $p:E \to B$ is called a \textbf{fibration} if, given any $k$-morphism $f:x\to y$ in $B$ and a lift $\tilde{x}$ of its source along $p$, there exists a lift $\tilde{f}:\tilde{x}\to\tilde{y}$ of $f$ along $p$. 
	
\end{definition}

We can rephrase this as saying that for every commutative square of the type formed by the solid arrows below, the dotted lift always exists, making both triangles commute.

$$\xymatrix{\theta^{(k-1)}\ar[r]^{\tilde{x}}\ar@{^{(}->}[d]_{s} & E\ar[d]^{p} \\ \theta^{(k)}\ar[r]_{f}\ar@{.>}[ru]^{\tilde{f}} & B}.$$ This amounts to saying that $p$ has the Right Lifting Property (RLP) with respect to the source inclusion map $s:\theta^{(k-1)}\hookrightarrow \theta^{(k)}$ for each $k=1, \cdots, n$. Here $\theta^{(k)}$ denotes the $n$-globular set represented by $k$.

\begin{remark}\label{folkfib}

A map of $n$-groupoids $p:E\to B$ is a fibration in the folk model structure on strict $n$-categories if every lifting problem of the form  \[\begin{tikzcd}
	{\theta^{(k)}} & X \\
	{I^{k+1}} & Y
	\arrow[hook', from=1-1, to=2-1]
	\arrow[from=1-1, to=1-2]
	\arrow[dotted, from=2-1, to=1-2]
	\arrow["f", from=1-2, to=2-2]
	\arrow[from=2-1, to=2-2]
\end{tikzcd}\] has a solution. Here $I^{k+1}$ can be thought of as a free walking coherent $(k+1)$-equivalence. It is defined in \cite{folk}, where it is denoted by $P^k$. Therefore, the question of whether the two notions of fibration agree is related to the possibility of extending equivalences to coherent equivalences. A mathoverflow question on this received no answers (\cite{question_fibrations}).

\begin{remark}
 
If the two notions of fibration agree, we can give a different proof of Theorem \ref{fibration} for all $n$, which we now sketch. In \cite{al20} it is shown that the category of strict $n$-categories equipped with the Crans-Gray tensor product and the folk model structure is a biclosed monoidal model category. This implies that the internal $\Hom$ functors $\Fun^{\lax}(-,\DD)$ and $\Fun^{\oplax}(-,\DD)$ send cofibrations to fibrations. From \cite{folk} one can deduce that an inclusion of presentations induces a cofibration between the presented $n$-categories. Therefore one can deduce that the restriction map on (op)lax functor categories is a folk fibration, from which it follows that the map of underlying $n$-groupoids is also a folk fibration. Since $\Map(\CC,\DD)$ is the underlying $n$-groupoid of $\Fun^{\oplax}(\CC,\DD)$, the Theorem would then follow.

However, once a theory of semistrict $4$-categories adpated to our string diagram calculus is in place, the proof of Theorem \ref{fibration} in the present paper will immeditely apply in that setting, which is not the case for the proof sketched above.
\end{remark}
 
\end{remark}

In the rest of this section we prove Theorem \ref{fibration}, splitting the proof into the various cases corresponding to the dimensions of the cells added to $\PP$ and the morphism to be lifted. We will repeatedly use the following colour code for distinguishing $k$-transfors in a mapping groupoid, which we already used in Section \ref{functorcat}: functors $F$ and $G$ are denoted in {\red} and {\blue}; natural transformations $\alpha$ and $\beta$ in {\green} and {\purple}; modifications $l$ and $m$ in {\orange} and {\lightblue}; perturbations $\A$ and $\B$ in {\olive} and {\teal}. The exception is when there is only one $k$-transfor involved, in which case we simply label it in black.

The proofs of the Lemmas in the rest of this section all follow the same line of reasoning, differing only in the explicit string diagram calculations. We have included all of them for completeness. It is probably useful to look at some to understand the general method, but then the reader can safely skip the rest.

\subsection{Adding a $0$-cell}

 \begin{lemma}
  
 Let $\PP$ be a presentation, $\CC$ a $4$-category and consider the presentation $\PP\cup\{Y\}$ obtained by adding a $0$-cell $Y$. Then the map $\Map(\PP\cup\{Y\},\CC)\to\Map(\PP,\CC)$ is a fibration of $4$-groupoids. 
  
 \end{lemma}
 
 \begin{proof}
  
 Suppose we have a $k$-morphism $x^k:a^{k-1}\to b^{k-1}$ in $\Map(\PP,\CC)$ together with an extension of $a$ to $\PP\cup\{Y\}$. We define $b_Y:=a_Y$ and $x_Y:=\id_{a_Y}$, which determines an extension of $x$ to $\PP\cup\{Y\}$.
 
 \end{proof}
 
 \subsection{Adding a $1$-cell}
 
 Let $\PP$ be a presentation, $\CC$ a $4$-category and consider the presentation $\PP\cup\{g\}$ obtained by adding a $1$-cell $g:X\to Y$, where $X$ and $Y$ are $0$-cells in $\PP$. We now prove that the map $\Map(\PP\cup\{g\},\CC)\to\Map(\PP,\CC)$ is a fibration of $4$-groupoids. We write $X=\includegraphics[scale=1.5,align=c]{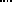}$, $Y=\includegraphics[scale=1.5,align=c]{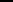}$ and $g=\includegraphics[scale=1.5,align=c]{morphisms/f.pdf}$. 
 
 \begin{lemma}
 	
 	 The map $\Map(\PP\cup\{g\},\CC)\to\Map(\PP,\CC)$ has the right lifting property with respect to $s:\theta^{(0)}\hookrightarrow \theta^{(1)}$.
 	
 \end{lemma}
 
 \begin{proof}
  
 Suppose we have a $1$-morphism $\alpha:F\to G$ in $\Map(\PP,\CC)$ together with an extension of $F$ to $\PP\cup\{g\}$. We want to define a $1$-morphism $G(g)=\includegraphics[scale=1,align=c]{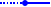}:\includegraphics[scale=1,align=c]{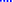}\to \includegraphics[scale=1,align=c]{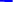}$ and a weakly invertible $2$-morphism $$\alpha_g:\includegraphics[scale=1,align=c]{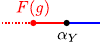}\to \includegraphics[scale=1,align=c]{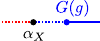}.$$ Since $\alpha_X$ is weakly invertible, the functor $$\includegraphics[scale=1,align=c]{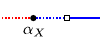}:\Hom\left(\includegraphics[scale=1,align=c]{fibration/1morph/GX.pdf},\includegraphics[scale=1,align=c]{fibration/1morph/GY.pdf}\right)\to\Hom\left(\includegraphics[scale=1,align=c]{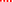},\includegraphics[scale=1,align=c]{fibration/1morph/GY.pdf}\right)$$ is essentially surjective, so we can pick a $1$-morphism $G(g)=\includegraphics[scale=1,align=c]{fibration/1morph/Gg.pdf}:\includegraphics[scale=1,align=c]{fibration/1morph/GX.pdf}\to \includegraphics[scale=1,align=c]{fibration/1morph/GY.pdf}$ such that $\includegraphics[scale=1,align=c]{fibration/1morph/a_g_s.pdf}\simeq \includegraphics[scale=1,align=c]{fibration/1morph/a_g_t.pdf}$ and define $\alpha_g$ to be any choice of such an equivalence.  
 
 \end{proof}
 
  \begin{lemma}
  	
 	 The map $\Map(\PP\cup\{g\},\CC)\to\Map(\PP,\CC)$ has the right lifting property with respect to $s:\theta^{(1)}\hookrightarrow \theta^{(2)}$.
  	
  \end{lemma}
  
  \begin{proof}
 
Consider two functors $F,G$ in $\Map(\PP\cup\{g\},\CC)$. Suppose we have also natural transformations $\alpha,\beta:F\to G$ and a modification $m:\alpha\to \beta$ in $\Map(\PP,\CC)$, together with an extension of $\alpha$ to $\PP\cup\{g\}$. We have weakly invertible $2$-morphisms \begin{center}\begin{tabular}{cc}$m_X= \includegraphics[scale=1.5,align=c]{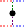}:\includegraphics[scale=1.5,align=c]{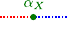}\to \includegraphics[scale=1.5,align=c]{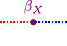}$ & $\alpha_g=  \includegraphics[scale=1.5,align=c]{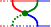}:\includegraphics[scale=1.5,align=c]{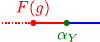}\to \includegraphics[scale=1.5,align=c]{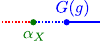}$ \\ $m_Y= \includegraphics[scale=1.5,align=c]{functorcat/2morph/m_y.pdf}: \includegraphics[scale=1.5,align=c]{functorcat/2morph/m_y_s.pdf}\to \includegraphics[scale=1.5,align=c]{functorcat/2morph/m_y_t.pdf}$ &  \end{tabular}\end{center} and we want to extend $\beta$ and $m$ to $\PP\cup\{g\}$, so we need a weakly invertible $2$-morphism \begin{center}$\beta_g=  \includegraphics[scale=1.5,align=c]{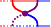}:\includegraphics[scale=1.5,align=c]{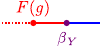}\to \includegraphics[scale=1.5,align=c]{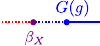}$\end{center} and a weakly invertible $3$-morphism \begin{center}$m_g:$ \includegraphics[scale=1.5,align=c]{functorcat/2morph/m_f_s.pdf} $\to$ \includegraphics[scale=1.5,align=c]{functorcat/2morph/m_f_t.pdf}.\end{center} The functor of $2$-categories $$\includegraphics[scale=1.4,align=c]{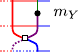}:\Hom\left(\includegraphics[scale=1.4,align=c]{fibration/2morph/b_g_s.pdf},\includegraphics[scale=1.4,align=c]{fibration/2morph/b_g_t.pdf}\right)\to\Hom\left(\includegraphics[scale=1.4,align=c]{fibration/2morph/a_g_s.pdf},\includegraphics[scale=1.4,align=c]{fibration/2morph/b_g_t.pdf}\right)$$ is given by composition with the weakly invertible $2$-morphism \includegraphics[scale=1]{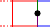}, and therefore it is essentially surjective. So we can pick a $2$-morphism $\beta_g:\beta_Y\circ F(g)\to G(g)\circ\beta_X$ such that \begin{center}\includegraphics[scale=1.5,align=c]{functorcat/2morph/m_f_s.pdf} $\simeq$ \includegraphics[scale=1.5,align=c]{functorcat/2morph/m_f_t.pdf}\end{center} and we pick $m_g$ to be any such equivalence.

\end{proof}

  \begin{lemma}
  	
 	 The map $\Map(\PP\cup\{g\},\CC)\to\Map(\PP,\CC)$ has the right lifting property with respect to $s:\theta^{(2)}\hookrightarrow \theta^{(3)}$.
  	
  \end{lemma}
 
 \begin{proof}
 Suppose we have a $3$-morphism $\A:l\to m$ in $\Map(\PP,\CC)$ together with an extension of $l$ to $\PP\cup\{g\}$. We want to define a $3$-morphism \begin{center}$m_g=$ \includegraphics[scale=1.5,align=c]{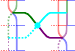} $:$ \includegraphics[scale=1.5,align=c]{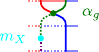} $\to$ \includegraphics[scale=1.5,align=c]{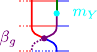}\end{center} and an invertible $4$-morphism \begin{center}$\A_g:$ \includegraphics[scale=1.5,align=c]{functorcat/3morph/a_f_s.pdf} $\to$ \includegraphics[scale=1.5,align=c]{functorcat/3morph/a_f_t.pdf}.\end{center} The functor $$\includegraphics[scale=1.3,align=c]{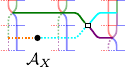}:\Hom\left(\includegraphics[scale=1.3,align=c]{fibration/3morph/m_f_s.pdf},\includegraphics[scale=1.3,align=c]{fibration/3morph/m_f_t.pdf}\right)\to\Hom\left(\includegraphics[scale=1.3,align=c]{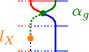},\includegraphics[scale=1.3,align=c]{fibration/3morph/m_f_t.pdf}\right)$$ is given by composition with a weakly invertible $3$-morphism, therefore it is essentially surjective. So we can pick a $3$-morphism \begin{center}$m_f=$ \includegraphics[scale=1.5,align=c]{fibration/3morph/m_f.pdf} $:$ \includegraphics[scale=1.5,align=c]{fibration/3morph/m_f_s.pdf} $\to$ \includegraphics[scale=1.5,align=c]{fibration/3morph/m_f_t.pdf}\end{center} such that \begin{center}\includegraphics[scale=1.5,align=c]{functorcat/3morph/a_f_s.pdf} $\simeq$ \includegraphics[scale=1.5,align=c]{functorcat/3morph/a_f_t.pdf}\end{center} and we can pick $\A_g$ to be any such isomorphism.
 
\end{proof}

  \begin{lemma}
  	
 	 The map $\Map(\PP\cup\{g\},\CC)\to\Map(\PP,\CC)$ has the right lifting property with respect to $s:\theta^{(3)}\hookrightarrow \theta^{(4)}$.
  	
  \end{lemma}
  
  \begin{proof}
 
Suppose we have a $4$-morphism $\Lambda:\A\to\B$ in $\Map(\PP,\CC)$ together with an extension of $\A$ to $\PP\cup\{g\}$. We want to define a $4$-morphism \begin{center}$\B_g=$ \includegraphics[scale=1.5,align=c]{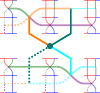} $:$ \includegraphics[scale=1.5,align=c]{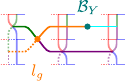} $\to$ \includegraphics[scale=1.5,align=c]{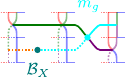}\end{center} such that \begin{center}\includegraphics[scale=1.5,align=c]{functorcat/4morph/f_f_lhs.pdf} $\xymatrix{\ar@{=}[r]^{\Lambda_g}&}$ \includegraphics[scale=1.5,align=c]{functorcat/4morph/f_f_rhs.pdf} .\end{center} The map $$\includegraphics[scale=1.2,align=c]{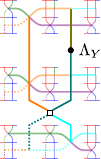}:\Hom\left(\includegraphics[scale=1.2,align=c]{fibration/4morph/b_f_s.pdf},\includegraphics[scale=1.2,align=c]{fibration/4morph/b_f_t.pdf}\right)\to\left(\includegraphics[scale=1.2,align=c]{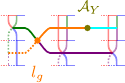},\includegraphics[scale=1.2,align=c]{fibration/4morph/b_f_t.pdf}\right)$$ is given by composition with an invertible $4$-morphism, therefore it is an isomorphism between sets of $4$-morphisms. In particular it is surjective, so we can pick a $4$-morphism $\B_g$ satisfying the desired relation.
 
 \end{proof}
 
 \subsection{Adding a $2$-cell}
 
 Let $\PP$ be a presentation, $\CC$ a $4$-category and consider the presentation $\PP\cup\{\zeta\}$ obtained by adding a $2$-cell $\zeta:f\to g$, where $f$ and $g$ are $1$-morphisms in $F_1(\PP)$. We will show that the map $\Map(\PP\cup\{\zeta\},\CC)\to\Map(\PP,\CC)$ is a fibration of $4$-groupoids. We denote $f$ by a dashed line and $g$ by a solid line and write $\zeta=$\includegraphics[scale=1.5,align=c]{morphisms/eta.pdf}.
 
 \begin{lemma}
  
The map $\Map(\PP\cup\{\zeta\},\CC)\to\Map(\PP,\CC)$ has the right lifting property with respect to $s:\theta^{(0)}\hookrightarrow \theta^{(1)}$.
  
 \end{lemma}
 
 \begin{proof}
  
 Suppose we have a $1$-morphism $\alpha:F\to G$ in $\Map(\PP,\CC)$ together with an extension of $F$ to $\PP\cup\{\zeta\}$. We want to define a $2$-morphism $G(\zeta):G(f)\to G(g)$ and a weakly invertible $3$-morphism \begin{center}$\alpha_{\zeta}:$   \includegraphics[scale=1.5,align=c]{functorcat/1morph/a_eta_s.pdf} $\to$ \includegraphics[scale=1.5,align=c]{functorcat/1morph/a_eta_t.pdf}.\end{center} The functor of $2$-categories \begin{center}\includegraphics[scale=1.5,align=c]{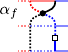} $:\Hom(G(f),G(g))\to\Hom(\alpha_Y\circ F(f),G(g)\circ\alpha_X)$\end{center} is given by first composing with the weakly invertible $1$-morphism $\alpha_X$ and then composing with the weakly invertible $2$-morphism $\alpha_f$. This implies this functor is essentially surjective, so we can choose $G(\zeta)\in\Hom(G(f),G(g))$ such that \begin{center}\includegraphics[scale=1.5,align=c]{functorcat/1morph/a_eta_s.pdf} $\simeq$ \includegraphics[scale=1.5,align=c]{functorcat/1morph/a_eta_t.pdf}\end{center} and we can pick $\alpha_{\zeta}$ to be any such equivalence. 
 
\end{proof}

 \begin{lemma}
 	
 	The map $\Map(\PP\cup\{\zeta\},\CC)\to\Map(\PP,\CC)$ has the right lifting property with respect to $s:\theta^{(1)}\hookrightarrow \theta^{(2)}$.
 	
 \end{lemma}
 
 \begin{proof}
 
Suppose we have a $2$-morphism $m:\alpha\to \beta$ in $\Map(\PP,\CC)$ together with an extension of $\alpha$ to $\PP\cup\{\zeta\}$. We want to define a $3$-morphism \begin{center}$\beta_{\zeta}=$ \includegraphics[scale=1.5,align=c]{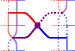} $:$  \includegraphics[scale=1.5,align=c]{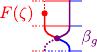} $\to$ \includegraphics[scale=1.5,align=c]{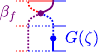}\end{center} and an invertible $4$-morphism \begin{center}$m_{\zeta}:$ \includegraphics[scale=1.5,align=c]{functorcat/2morph/m_eta_s.pdf} $\to$ \includegraphics[scale=1.5,align=c]{functorcat/2morph/m_eta_t.pdf}.\end{center} The functor $$\includegraphics[scale=1.25,align=c]{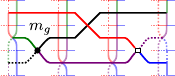} :\Hom\left(\includegraphics[scale=1.25,align=c]{fibration/2morph/b_eta_s.pdf}, \includegraphics[scale=1.25,align=c]{fibration/2morph/b_eta_t.pdf}\right)\to\Hom\left(\includegraphics[scale=1.25,align=c]{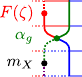},\includegraphics[scale=1.25,align=c]{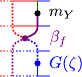} \right)$$ is given by first composing with a weakly invertible $2$-morphism and then with a weakly invertible $3$-morphism, therefore it is essentially surjective. So we can pick a $3$-morphism \begin{center}$\beta_{\zeta}:$ \includegraphics[scale=1.5,align=c]{fibration/2morph/b_eta_s.pdf} $\to$ \includegraphics[scale=1.5,align=c]{fibration/2morph/b_eta_t.pdf}\end{center} such that \begin{center}\includegraphics[scale=1.5,align=c]{functorcat/2morph/m_eta_s.pdf} $\simeq$ \includegraphics[scale=1.5,align=c]{functorcat/2morph/m_eta_t.pdf}\end{center} and we can pick $m_{\zeta}$ to be any such isomorphism.

\end{proof}

\begin{lemma}
	
	The map $\Map(\PP\cup\{\zeta\},\CC)\to\Map(\PP,\CC)$ has the right lifting property with respect to $s:\theta^{(2)}\hookrightarrow \theta^{(3)}$.
	
\end{lemma}

\begin{proof}
 
 Suppose we have a $3$-morphism $\A:l\to m$ in $\Map(\PP,\CC)$ together with an extension of $l$ to $\PP\cup\{\zeta\}$. We want to define a $4$-morphism \begin{center}$m_{\zeta}=$ \includegraphics[scale=1.5,align=c]{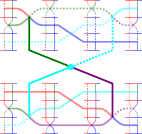} $:$\includegraphics[scale=1.5,align=c]{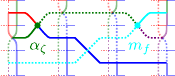} $\to$ \includegraphics[scale=1.5,align=c]{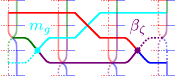}\end{center} such that \begin{center}\includegraphics[scale=1.5,align=c]{functorcat/3morph/a_eta_lhs.pdf}$\xymatrix{\ar@{=}[r]^{\A_{\zeta}}&}$\includegraphics[scale=1.5,align=c]{functorcat/3morph/a_eta_rhs.pdf}.\end{center} The map $$\Hom\left(\includegraphics[scale=.85,align=c]{fibration/3morph/m_eta_s.pdf},\includegraphics[scale=.85,align=c]{fibration/3morph/m_eta_t.pdf}\right)\to\Hom\left(\includegraphics[scale=0.85,align=c]{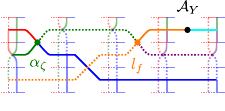},\includegraphics[scale=0.85,align=c]{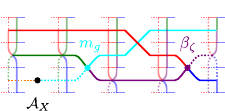}\right)$$ given by $$\includegraphics[scale=1.5,align=c]{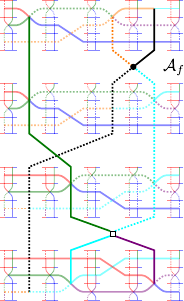}$$ can be obtained by composing with weakly invertible mophisms, therefore it is surjective. So we can pick $m_{\zeta}$ satisfying the desired relation.
 
\end{proof}

\begin{lemma}
	
	The map $\Map(\PP\cup\{\zeta\},\CC)\to\Map(\PP,\CC)$ has the right lifting property with respect to $s:\theta^{(3)}\hookrightarrow \theta^{(4)}$.
	
\end{lemma}

\begin{proof}
 
 Suppose we have a $4$-morphism $\Lambda:\A\to \B$ in $\Map(\PP,\CC)$ together with an extension of $\A$ to $\PP\cup\{\zeta\}$, which means we have the relation  \begin{center}\includegraphics[scale=1.5,align=c]{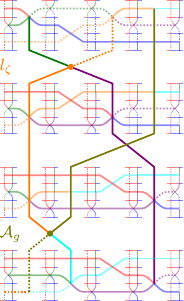} $\xymatrix{\ar@{=}[r]^{\A_{\zeta}}&}$ \includegraphics[scale=1.5,align=c]{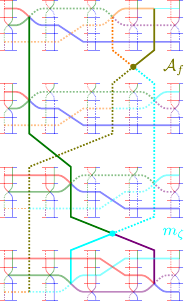} .\end{center} We want to extend $\B$ to $\PP\cup\{\zeta\}$, so we need to prove that $\B$ satisfies the relation \begin{center}\includegraphics[scale=1.5,align=c]{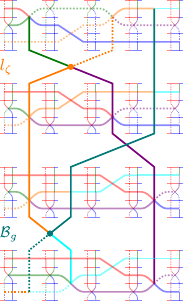} $\xymatrix{\ar@{=}[r]^{\B_{\zeta}}&}$ \includegraphics[scale=1.5,align=c]{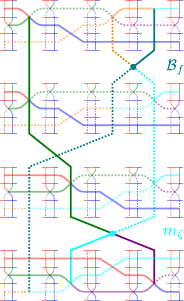} .\end{center} Using the relations $\Lambda_f$, $\Lambda_g$ and $\A_{\zeta}$ we have the following proof, where we use $\stackrel{h}{=}$ when two $4$-diagrams are related by homotopy geneators.
 
  \begin{longtable}{llllll}
   
 \includegraphics[scale=1,align=c]{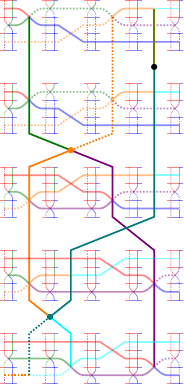}  & $\xymatrix{\ar@{=}[r]^{h}&}$ & \includegraphics[scale=1,align=c]{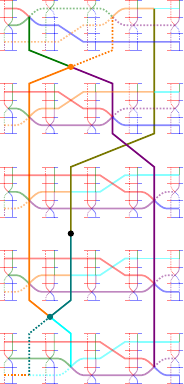} & $\xymatrix{\ar@{=}[r]^{\Lambda_g^{-1}}&}$ & \includegraphics[scale=1,align=c]{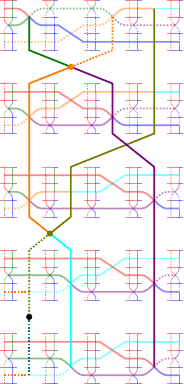} & $\xymatrix{\ar@{=}[r]^{\A_{\zeta}}&}$ 
 
 \\
 
 \\
 
  \includegraphics[scale=1,align=c]{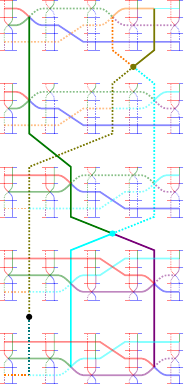} & $\xymatrix{\ar@{=}[r]^{h}&}$ &  \includegraphics[scale=1,align=c]{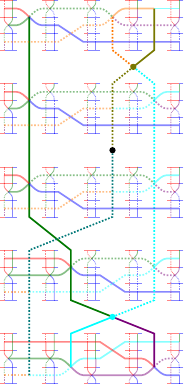} & $\xymatrix{\ar@{=}[r]^{\Lambda_f}&}$ &  \includegraphics[scale=1,align=c]{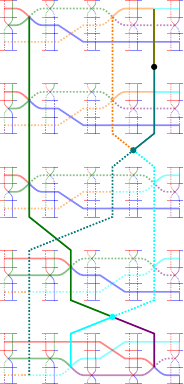} & $.$

  \end{longtable}
  
The map $$\Hom\left(\includegraphics[scale=0.76,align=c]{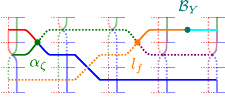},\includegraphics[scale=0.76,align=c]{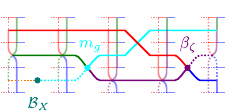}\right)\to\Hom\left(\includegraphics[scale=0.76,align=c]{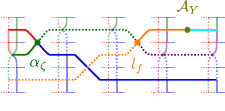},\includegraphics[scale=0.76,align=c]{fibration/4morph/b_eta_fill_3.pdf}\right)$$ given by $$\includegraphics[scale=1.5,align=c]{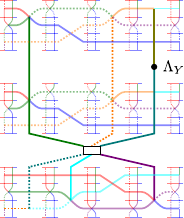}$$ can be obtained by composition with an invertible $4$-morphism, therefore it is an isomorphism between sets of $4$-morphisms. In particular it is injective, so the desired relation follows.
   
 \end{proof}
 
 \subsection{Adding a $3$-cell}
 
 Let $\PP$ be a presentation, $\CC$ a $4$-category and consider the presentation $\PP\cup\{t\}$ obtained by adding a $3$-cell $t:\eta\to \zeta$, where $\eta$ and $\zeta$ are $2$-morphisms in $F_2(\PP)$. We will show that the map $\Map(\PP\cup\{t\},\CC)\to\Map(\PP,\CC)$ is a fibration of $4$-groupoids. We denote $\eta$ by a dashed line and $\zeta$ by a solid line and write $t=$\includegraphics[scale=1.5,align=c]{morphisms/s.pdf}.

\begin{lemma}
	
	The map $\Map(\PP\cup\{t\},\CC)\to\Map(\PP,\CC)$ has the right lifting property with respect to $s:\theta^{(0)}\hookrightarrow \theta^{(1)}$.
	
\end{lemma}
  
 \begin{proof}
  
Suppose we have a $1$-morphism $\alpha:F\to G$ in $\Map(\PP,\CC)$ together with an extension of $F$ to $\PP\cup\{t\}$. We want to define a $3$-morphism $G(t):G(\eta)\to G(\zeta)$ and an invertible $4$-morphism \begin{center}$\alpha_{t}:$ \includegraphics[scale=1.5,align=c]{functorcat/1morph/a_t_s.pdf} $\to$ \includegraphics[scale=1.5,align=c]{functorcat/1morph/a_t_t.pdf} .\end{center} The functor \begin{center}\includegraphics[scale=1.5,align=c]{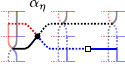} $:\Hom(G(\eta),G(\zeta))\to\Hom\left(\includegraphics[scale=1.5,align=c]{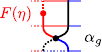},\includegraphics[scale=1.5,align=c]{functorcat/1morph/a_eta_t.pdf}\right)$\end{center} is given by composing with weakly invertible morphisms. This implies that it is essentially surjective, so we can choose $G(t):G(\eta)\to G(\zeta)$ such that \begin{center}\includegraphics[scale=1.5,align=c]{functorcat/1morph/a_t_s.pdf} $\simeq$ \includegraphics[scale=1.5,align=c]{functorcat/1morph/a_t_t.pdf}\end{center} and we can pick $\alpha_t$ to be any such isomorphism.

\end{proof}

\begin{lemma}
	
	The map $\Map(\PP\cup\{t\},\CC)\to\Map(\PP,\CC)$ has the right lifting property with respect to $s:\theta^{(1)}\hookrightarrow \theta^{(2)}$.
	
\end{lemma}

\begin{proof}
 
Suppose we have a $2$-morphism $m:\alpha\to \beta$ in $\Map(\PP,\CC)$ together with an extension of $\alpha$ to $\PP\cup\{t\}$. We want to define a $4$-morphism \begin{center}$\beta_{t}=$ \includegraphics[scale=1.5,align=c]{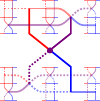} $:$  \includegraphics[scale=1.5,align=c]{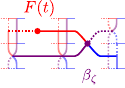} $\to$ \includegraphics[scale=1.5,align=c]{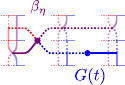}\end{center} such that \begin{center}\includegraphics[scale=1.5,align=c]{functorcat/2morph/m_t_lhs.pdf} $\xymatrix{\ar@{=}[r]^{m_t}&}$ \includegraphics[scale=1.5,align=c]{functorcat/2morph/m_t_rhs.pdf} .\end{center} The map $$\Hom\left(\includegraphics[scale=0.98,align=c]{fibration/2morph/b_t_s.pdf},\includegraphics[scale=.98,align=c]{fibration/2morph/b_t_t.pdf}\right)\to\Hom\left(\includegraphics[scale=.98,align=c]{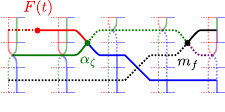},\includegraphics[scale=.98,align=c]{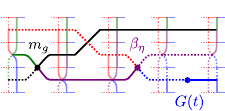}\right)$$ given by $$\includegraphics[scale=1.5,align=c]{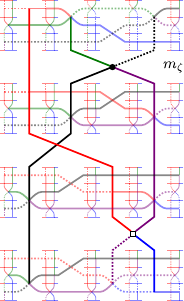}$$ can be obtained by composition with weakly invertible morphisms, so it is surjective. This implies we can choose $\beta_t$ so that the desired equality holds.
  
\end{proof}

\begin{lemma}
	
	The map $\Map(\PP\cup\{t\},\CC)\to\Map(\PP,\CC)$ has the right lifting property with respect to $s:\theta^{(2)}\hookrightarrow \theta^{(3)}$.
	
\end{lemma}

\begin{proof}  
  
Suppose we have a $3$-morphism $\A:l\to m$ in $\Map(\PP,\CC)$ together with an extension of $l$ to $\PP\cup\{t\}$, which means we have the relation  \begin{center}\includegraphics[scale=1.5,align=c]{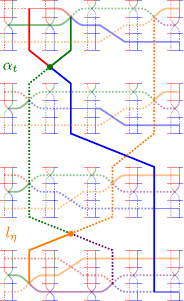} $\xymatrix{\ar@{=}[r]^{l_t}&}$ \includegraphics[scale=1.5,align=c]{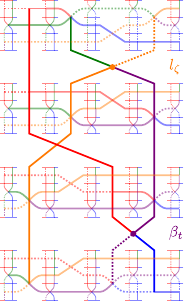} .\end{center} We want to extend $m$ to $\PP\cup\{t\}$, so we need to prove that $m$ satisfies the relation \begin{center}\includegraphics[scale=1.5,align=c]{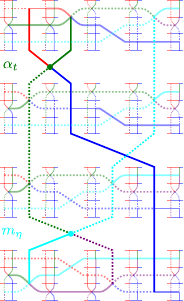} $\xymatrix{\ar@{=}[r]^{m_t}&}$ \includegraphics[scale=1.5,align=c]{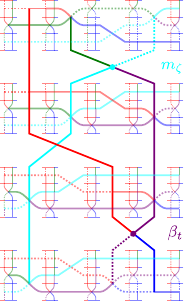} .\end{center} Using the relations $\A_{\eta}$, $\A_{\zeta}$ and $l_t$ we get

  \begin{longtable}{llllll}
   
 \includegraphics[scale=1,align=c]{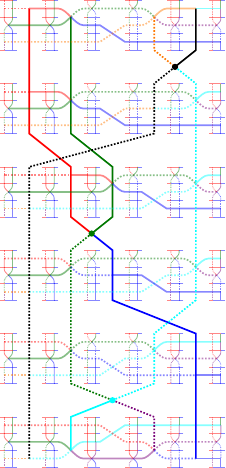}  & $\stackrel{h}{=}$ & \includegraphics[scale=1,align=c]{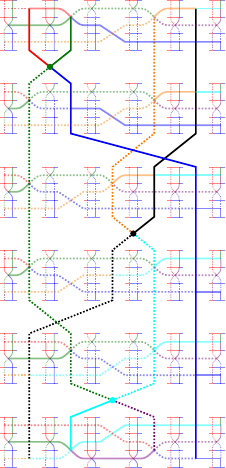} & $\stackrel{\A_{\eta}^{-1}}{=}$ &  \includegraphics[scale=1,align=c]{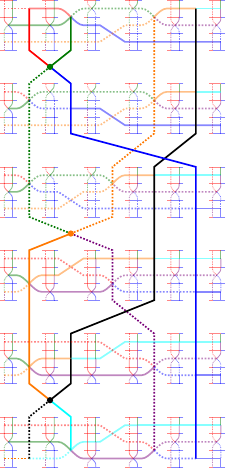} & $\stackrel{h}{=}$ 
 
\\

\\

\includegraphics[scale=1,align=c]{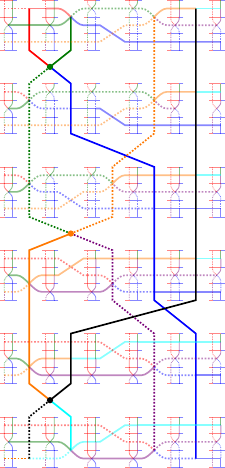} & $\stackrel{l_t}{=}$ & \includegraphics[scale=1,align=c]{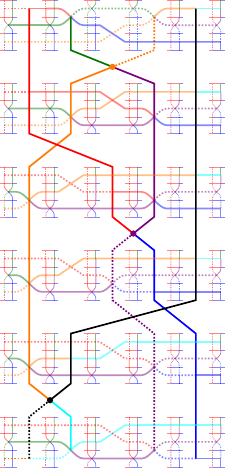} & $\stackrel{h}{=}$ &  \includegraphics[scale=1,align=c]{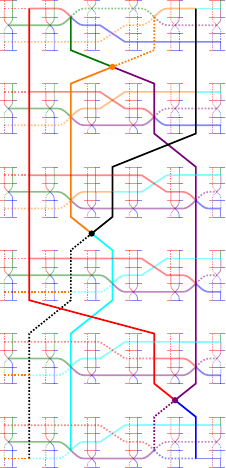} & $\stackrel{\A_{\zeta}}{=}$
 
 \\
 
 \\
 
 \includegraphics[scale=1,align=c]{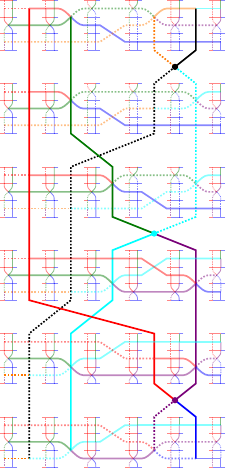} & $\stackrel{h}{=}$ &  \includegraphics[scale=1,align=c]{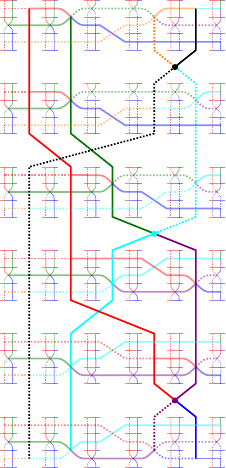} & $.$ & &

  \end{longtable}
  
The map $$\Hom\left(\includegraphics[scale=0.69,align=c]{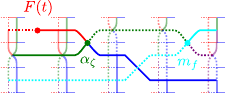},\includegraphics[scale=0.69,align=c]{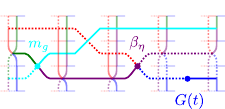}\right)\to\Hom\left(\includegraphics[scale=0.69,align=c]{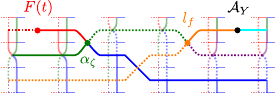},\includegraphics[scale=0.69,align=c]{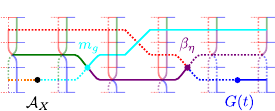}\right)$$ given by $$\includegraphics[scale=1.5,align=c]{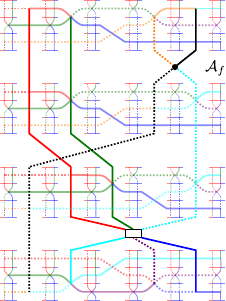}$$ can be obtained by composition with weakly invertible morphisms, so it is an isomorphism between sets of $4$-morphisms. In particular it is injective, so the desired relation follows.
  
 \end{proof}
 
 \subsection{Adding a $4$-cell}
  
 Let $\PP$ be a presentation, $\CC$ a $4$-category and consider the presentation $\PP\cup\{\W\}$ obtained by adding a $4$-cell $\W:s \to t$, where $s$ and $t$ are $3$-morphisms in $F_3(\PP)$. We will show that the map $\Map(\PP\cup\{\W\},\CC)\to\Map(\PP,\CC)$ is a fibration of $4$-groupoids. We denote $s$ by a dashed line and $t$ by a solid line and write $$\W=\includegraphics[scale=1.5,align=c]{morphisms/w.pdf}.$$

\begin{lemma}
	
	The map $\Map(\PP\cup\{\W\},\CC)\to\Map(\PP,\CC)$ has the right lifting property with respect to $s:\theta^{(0)}\hookrightarrow \theta^{(1)}$.
	
\end{lemma}
 
 \begin{proof}
  
 Suppose we have a $1$-morphism $\alpha:F\to G$ in $\Map(\PP,\CC)$ together with an extension of $F$ to $\PP\cup\{\W\}$. We want to define a $4$-morphism $G(\W):G(s)\to G(t)$ such that \begin{center} \includegraphics[scale=1.5,align=c]{functorcat/1morph/a_a_lhs.pdf} $\xymatrix{\ar@{=}[r]^{\alpha_{\W}}&}$  \includegraphics[scale=1.5,align=c]{functorcat/1morph/a_a_rhs.pdf} .\end{center} The map $$\includegraphics[scale=1.5,align=c]{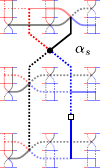}:\Hom(G(s),G(t))\to\Hom\left(\includegraphics[scale=1.5,align=c]{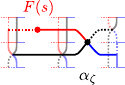},\includegraphics[scale=1.5,align=c]{functorcat/1morph/a_t_t.pdf}\right)$$ is given by composing with weakly invertible morphisms. This implies that it is surjective, so we can choose $G(\W):G(s)\to G(t)$ such that \begin{center}\includegraphics[scale=1.5,align=c]{functorcat/1morph/a_a_lhs.pdf} $\xymatrix{\ar@{=}[r]^{\alpha_{\W}}&}$ \includegraphics[scale=1.5,align=c]{functorcat/1morph/a_a_rhs.pdf} .\end{center}
 
 \end{proof}
 
 \begin{lemma}
 	
 	The map $\Map(\PP\cup\{\W\},\CC)\to\Map(\PP,\CC)$ has the right lifting property with respect to $s:\theta^{(1)}\hookrightarrow \theta^{(2)}$.
 	
 \end{lemma}
 
 \begin{proof}
 
 Suppose we have a $2$-morphism $m:\alpha\to \beta$ in $\Map(\PP,\CC)$ together with an extension of $\alpha$ to $\PP\cup\{\W\}$, which means we have \begin{center}\includegraphics[scale=1.5,align=c]{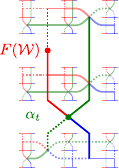}$\xymatrix{\ar@{=}[r]^{\alpha_{\W}}&}$\includegraphics[scale=1.5,align=c]{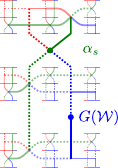}.\end{center} We want to extend $\beta$ to $\W$, so we want to show that \begin{center}\includegraphics[scale=1.5,align=c]{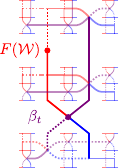}$\xymatrix{\ar@{=}[r]^{\beta_{\W}}&}$\includegraphics[scale=1.5,align=c]{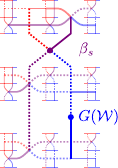}.\end{center}
 
 Using the relations $m_s$, $m_t$ and $\alpha_{\W}$, we get 
 
  \begin{longtable}{llllll}
   
 \includegraphics[scale=1,align=c]{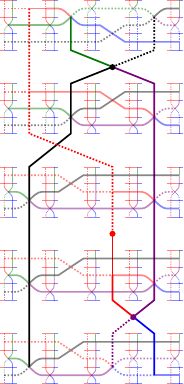}  & $\xymatrix{\ar@{=}[r]^{h}&}$ & \includegraphics[scale=1,align=c]{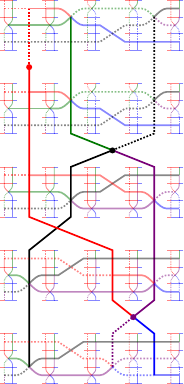} & $\xymatrix{\ar@{=}[r]^{m_{t}^{-1}}&}$ & \includegraphics[scale=1,align=c]{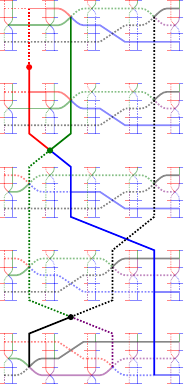} & $\xymatrix{\ar@{=}[r]^{\alpha_{\W}}&}$ 
 
 \\
 
 \\
 
  \includegraphics[scale=1,align=c]{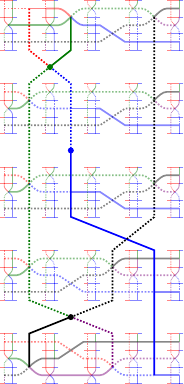} & $\xymatrix{\ar@{=}[r]^{h}&}$ &  \includegraphics[scale=1,align=c]{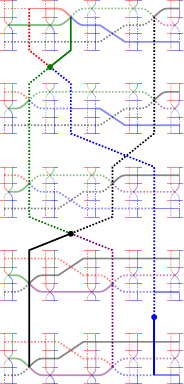} & $\xymatrix{\ar@{=}[r]^{m_{s}}&}$ &  \includegraphics[scale=1,align=c]{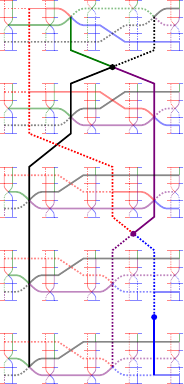} & $.$

  \end{longtable}
  
 The map $$\Hom\left(\includegraphics[scale=0.9,align=c]{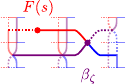},\includegraphics[scale=0.9,align=c]{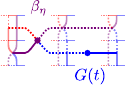}\right)\to\Hom\left(\includegraphics[scale=0.9,align=c]{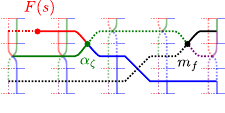},\includegraphics[scale=0.9,align=c]{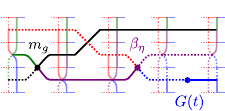}\right)$$ given by $$\includegraphics[scale=1.5,align=c]{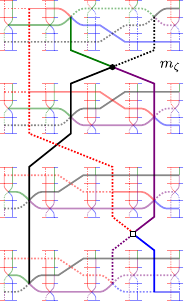}$$ can be obtained by composing with weakly invertible morphisms, so it is an isomorphism and we get the desired relation.
  
 \end{proof}
 
 \subsection{Adding a $5$-cell}
  
 Let $\PP$ be a presentation, $\CC$ a $4$-category and consider the presentation $\PP\cup\{\Phi\}$ obtained by adding a $5$-cell $\Phi:\VV \to \W$, where $\VV$ and $\W$ are $4$-morphisms in $F_4(\PP)$. We will show that the map $\Map(\PP\cup\{\Phi\},\CC)\to\Map(\PP,\CC)$ is a fibration of $4$-groupoids. 
  
  \begin{lemma}
  	
  	The map $\Map(\PP\cup\{\Phi\},\CC)\to\Map(\PP,\CC)$ has the right lifting property with respect to $s:\theta^{(0)}\hookrightarrow \theta^{(1)}$.
  	
  \end{lemma}
 
 \begin{proof}
  
 Suppose we have a $1$-morphism $\alpha:F\to G$ in $\Map(\PP,\CC)$ together with an extension of $F$ to $\PP\cup\{\Phi\}$. Recall that in a $4$-category a $5$-cell corresponds to a relation, so we have $F(\VV)=F(\W)$. We want to extend $G$ to $\PP\cup\{\Phi\}$, so we want to show that $G(\VV)=G(\W)$. Using the relations $\alpha_{\VV}$ and $\alpha_{\W}$ we get

  \begin{longtable}{lllllll}
   
 \includegraphics[scale=1,align=c]{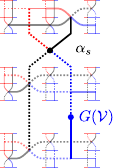}  & $\xymatrix{\ar@{=}[r]^{\alpha_{\VV}^{-1}}&}$ & \includegraphics[scale=1,align=c]{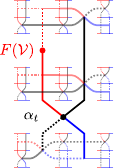} & $\xymatrix{\ar@{=}[r]^{F(\Phi)}&}$ & \includegraphics[scale=1,align=c]{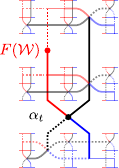} & $\xymatrix{\ar@{=}[r]^{\alpha_{\W}}&}$ & \includegraphics[scale=1,align=c]{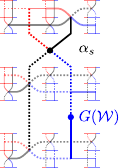}.

  \end{longtable}

 The map $$\includegraphics[scale=1.5,align=c]{fibration/1morph/a_a_rhs.pdf}:\Hom(G(s),G(t))\to\Hom\left(\includegraphics[scale=1.5,align=c]{fibration/1morph/a_s_s.pdf},\includegraphics[scale=1.5,align=c]{functorcat/1morph/a_t_t.pdf}\right)$$ is invertible, as it is given by composing with weakly invertible morphisms. This implies that it is injective, so $G(\VV)=G(\W)$.
  
 \end{proof} 
 
This completes the proof of Theorem \ref{fibration} in the case where $\QQ$ is obtained from $\PP$ by adding a cell. Since a composite of fibrations is a fibration, the general case follows.

\end{document}